\documentclass[11pt]{article}
\usepackage{amsmath}
\usepackage{algorithm}
\usepackage{algpseudocode}

\usepackage[english]{babel}

\usepackage[letterpaper,top=2cm,bottom=2cm,left=3cm,right=3cm,marginparwidth=1.75cm]{geometry}

\usepackage[utf8]{inputenc}
\usepackage{amsmath}
\usepackage{amssymb}
\usepackage{amsthm}
\usepackage{multicol}
\usepackage{bm}
\usepackage{multirow}
\usepackage{amssymb}
\usepackage{algpseudocode} 
\usepackage{algorithm}
\usepackage[english]{babel}
\usepackage{adjustbox}
\usepackage[T1]{fontenc}
\usepackage{url}            
\usepackage{booktabs}       
\usepackage{amsfonts}       
\usepackage{nicefrac}       
\usepackage{microtype}      
\usepackage{tikz}
\usetikzlibrary{positioning}
\usepackage{graphicx}
\usepackage[nottoc,notlot,notlof]{tocbibind}
\usepackage{caption}
\usepackage{subcaption}
\usepackage{mathabx}
\usepackage[inline]{enumitem}
\usepackage{mathtools}
\usepackage{rotating}
\usepackage{csquotes}

\newtheorem{theorem}{Theorem}
\newtheorem{assumption}{Assumption}
\newtheorem{proposition}{Proposition}
\newtheorem{corollary}{Corollary}
\newtheorem{lemma}{Lemma}
\newtheorem{remark}{Remark}

\def\Pruno{\mathsf{P}_1}
\def\Prdue{\mathsf{P}_2}
\def\Csf{\mathsf{C}}
\def\Dsf{\mathsf{D}}
\def\Bsf{\mathsf{B}}
\def\Xsf{\mathsf{X}}
\def\Ceta{\Csf_\eta}
\def\Cu{\Csf_u}
\def\Cw{\Csf_w}
\def\Lift{\mathsf{L}}

\def\XC{\mathcal{X}}
\def\YC{\mathcal{Y}}
\def\KC{\mathcal{K}}
\def\AC{\mathcal{A}}
\def\BC{\mathcal{B}}

\def\setKR{\mathcal{P}_R}

\def\DivBp{\mathcal{B}_p}

\def \erre{\mathbb{R}}
\def \ER{\mathcal{E}_R}
\def \HR{\ER^c}

\def \GR{G_R}

\def\Up{\psi_p}
\def\Udue{\psi_2}
\def\Upo{\psi^\prime_{p,0}}
\def\Upi{\psi^\prime_{p,\infty}}
\def \Isubl{\mathcal C_0}

\def \FeasR{\mathcal{F}_R}
\def \FeasRp{\mathcal{F}_{p,R}}

\def\Tij{T_{ij}}
\def\Tijp{\Tij^{(p)}}

\def\minr{c_{i,\min}}
\def\minc{c_{\min,j}}
\def\Iphi{I_{\phi}}

\def\ctrun{c^{(R)}}

\def\dom{\mathrm{dom}}
\def\supp{\mathrm{supp}}
\DeclareMathOperator{\argmin}{argmin}
\DeclareMathOperator{\ri}{ri}
\newcommand{\pair}[2]{\langle #1,#2\rangle}

\def\phipz{\phi^\prime_0}
\def\phipi{\phi_\infty^\prime}

\def\one{\mathbf 1}
\newcommand{\Mplus}{\mathcal M_+}
\def\ind{\iota}

\def \UTp{\mathrm{UT}_{p,c}}

\def \CC{\mathcal{C}}

\def \UTphi {\mathrm{UT}_{\phi,c}}
\def \UTRphint {\mathrm{UT}_{\phi,\ctrun}}

\def \Dp{\mathrm{D}_p}
\def \Dphi{\mathrm{D}_\phi}

\usepackage{mathrsfs}
\numberwithin{equation}{section}

\providecommand{\scalarp}[1]{\langle#1\rangle}

\newcommand{\R}{\ensuremath \mathbb{R}}

\newcommand{\ZC}{\ensuremath \mathcal{Z}}

\title{Transport Less, Solve Faster: Exploiting sparsity in algorithm design for Unbalanced Optimal Transport}

\author{%
  Gennaro Auricchio\thanks{%
    Department of Mathematics ``Tullio Levi-Civita'',
    University of Padova, Padova, Italy.
    Emails: \texttt{gennaro.auricchio@unipd.it},
    \texttt{rinaldi@math.unipd.it},
    \texttt{giovanni.spisso@phd.unipd.it}.%
  }
  \and
  Francesco Rinaldi\footnotemark[1]
  \and
  Saverio Salzo\thanks{%
    Department of Computer, Control and Management Engineering
    ``Antonio Ruberti'',
    Sapienza University of Rome, Rome, Italy.
    Email: \texttt{saverio.salzo@uniroma1.it}.%
  }
  \and
  Giovanni Spisso\footnotemark[1]
}

\date{}

\usepackage[backend=biber,style=numeric, giveninits=true, url=false, isbn=false, maxbibnames=99]{biblatex}
\AtEveryBibitem{%
  \clearname{editor}%
  \clearname{editora}%
  \clearname{editorb}%
}

\graphicspath{{figures/}}

\begin{document}

\maketitle

\begin{abstract}
    In this paper, we study the computational aspects of the Unbalanced Optimal Transport (UOT) problem between two given positive measures $\mu$ and $\nu$.
    The UOT is a variant to the classic Optimal Transportation problem in which the linear constraints modeling mass preservation are relaxed and replaced by a convex penalty $\Dphi$ in the objective function.
    First, we show that any UOT problem admits a discrete solution as long as $\mu$ and $\nu$ are discrete, allowing us to express the optimality conditions characterizing the solution to the UOT problem in an easy-to-handle form.
    By exploiting the optimality conditions we provide an upper and lower bound on the marginals of the solution as well as a bound on the mass that any solution assigns to any couple of points.
    The bounds on the mass assigned to each couple of points reveal that, under suitable conditions, the support of any solution to the UOT is contained within a compact set that does not depend on the measures $\mu$ and $\nu$, but only on the penalty function $\Dphi$ and the cost function $c$.
    We then introduce a truncated version of the UOT problem that fully leverages the sparsity of the optimal solution by introducing a parameter $R$ that saturates $c$.
    We show that any solution to the truncated problem naturally yields an approximation of the optimal solution for the standard UOT problem and bound the primal sub-optimality gap between the real solution and the solution induced by the truncated problem via a decreasing power-like function of $R$.
    Lastly, we study the convergence speed of the Pairwise Frank-Wolfe algorithm when applied to the truncated UOT problem and show that such routine is able to identify the support of the solution within a finite number of iterations.
\end{abstract}

\section{Introduction}

Optimal Transport (OT) is a standard way to lift any distance defined over a Polish space $(\Omega,d)$ to the space of probability measures defined over $\Omega$ \cite{villani2008optimal,santambrogio2015optimal}.
Due to their geometrical properties the distances induced by OT problems, i.e., the Wasserstein distances, have been proven to be a valuable tool to address several key sub-tasks in central computer science problems, e.g., Machine Learning \cite{scagliotti2023normalizing,Cuturi2014,Frogner2015,pmlr-v70-arjovsky17a,pereira2025survey}, Computer Vision \cite{Rubner1998,rubner2000earth,Pele2009,Granger2024b}, and Natural Language Processing \cite{kusner2015word,rolet2016fast}.
Owing to the growing interest for OT-related distances, the demand for fast and efficient computation methods grew accordingly.
Indeed, in its modern formulation, Wasserstein Distance induced by a cost function $c:\XC\times\YC\to [0,+\infty[$ between two probability measures $\mu\in\mathcal{P}(\XC)$ and $\nu\in\mathcal{P}(\YC)$ is defined as the minimum of the following linear problem
\begin{equation}
\label{eq:kantorovich}
\min_{\pi \in \Pi(\mu,\nu)}
\int_{\mathcal{X}\times\mathcal{Y}} c(x,y)\, \pi(dx,dy),
\end{equation}
where $\Pi(\mu,\nu) = \left\{ \pi \in \mathcal{P}(\mathcal{X}\times\mathcal{Y}) \,\middle|\, \pi^{(1)}=\mu,\; \pi^{(2)}=\nu \right\}$ with $\pi^{(1)}$ and $\pi^{(2)}$ being the marginal of $\pi$ over $\XC$ and the marginal of $\pi$ over $\YC$, respectively.
For discrete probability measures, problem \eqref{eq:kantorovich} becomes a Linear Programming (LP) problem, whose computational cost becomes prohibitive when the dimension of the problem increases \cite{kantorovich2006translocation}.
The demand for fast algorithms capable of solving the OT problem has spurred extensive research aimed at reducing the computational cost of these methods by considering a Linear Programming problem that exploits the structure of the problem \cite{LingOkada2007,Auricchio2018} or by approaching an entropic regularization version of the problem \cite{Cuturi2013,Solomon2015,Altschuler2017}.
Despite the success of the Wasserstein Distances, they are limited as they compare only positive measures with the same total mass.
This factor represents a flaw in any application where the total mass of the object under study is important, e.g., in shape analysis \cite{chizat2018interpolating,Solomon2015} or in persistence diagrams comparison \cite{vidal2019progressive,lacombe2018large}.
Moreover, since the OT problem imposes the conservation of the mass, the resulting Wasserstein Distances are sensitive to noise \cite{sato2020fast}, hindering the applicability of the Wasserstein Distances to fields such as robust learning \cite{balaji2020robust}, domain adaptation \cite{fatras2021unbalanced}, imaging \cite{lellmann2014imaging}, biology \cite{schiebinger2019optimal}, and text matching \cite{tang2023learning}.
To overcome these shortcomings, there have been several proposals to extend the Wasserstein Distances to the set of positive measures, such as partial optimal transport \cite{figalli2010optimal}, the Wasserstein Fisher-Rao distance \cite{hanin1992kantorovich,chizat2018interpolating}, or  the Kantorovich Rubinstein distance \cite{lellmann2014imaging}.
In this paper, we focus on the proposal presented by Liero, Mielke and Savaré, where the hard constraint on the marginal of the transportation plan between the two probability measures is replaced by a suitable Csiszár penalty term in the objective function \cite{liero2018optimal}.
More formally, given two positive measures $\mu$ and $\nu$, the problem we consider is as follows
\begin{equation*}
\min_{\pi\in \mathcal{M}_+(\XC\times \YC)} \UTphi(\pi) := \min_{\pi\in \mathcal{M}_+(\XC\times \YC)} \int_{\XC\times\YC} c(x,y) \pi(dx, dy) + \Dphi(\pi^{(1)} \,\vert\, \mu) + \Dphi(\pi^{(2)} \,\vert\, \nu),
\end{equation*}
where $\Dphi$ is the Csiszár divergence induced by a convex function $\phi$.
Unlike the classic Wasserstein Distance, the unbalanced optimal transportation problem proposed in \cite{liero2018optimal} can be naturally formulated as the minimization of a convex function over a convex set, enabling the use of classic routines for convex programming.

\subsection{Related Works}

As per the standard OT problem, most approaches for the UOT rely on a convex relaxation of the objective function that is then handled through the Sinkhorn algorithm \cite{Cuturi2013,sinkhorn1967concerning}.
The first Sinkhorn-like algorithm for unbalanced optimal transport was proposed in \cite{chizat2018scaling}, where the authors showed how to extend the entropic regularization techniques used for the classic OT problem. 
This approach was subsequently improved in \cite{schmitzer2019stabilized}, where the authors enhanced the numerical stability of the method.
For the one dimensional case, the state-of-the-art solver of the UOT was introduced in \cite{sejourne2022faster}, where the authors considered a translation-invariant version of the dual UOT problem that can be efficiently solved using the Frank-Wolfe algorithm.
The results presented in \cite{sejourne2022faster} have led to the study of a Sliced Unbalanced Wasserstein Distances \cite{bonet2024slicing} in which the UOT distance is computed by taking the average over all the UOT distances computed over all the possible one-dimensional projections. 
A different approach based on entropic regularization is presented in \cite{scetbon2023unbalanced}, where the authors impose an additional low-rank constraint to the UOT problem, enabling a more compact representation of the feasible set of transport plans.
An alternative regularization approach was presented in \cite{blondel2018smooth,nguyen2023unbalanced}, where the authors define a UOT problem regularization via an $\ell_2$ penalty on the transport plan rather than an entropic one, yielding a strongly convex objective function that can be tackled via standard gradient schemes \cite{nguyen2023unbalanced}.
In \cite{chapel2021unbalanced}, the authors approach the solution of the UOT problem by recasting it as a convex penalized linear regression problem, which has been proven to be efficient when handling specific Csiszár divergences.
Lastly, in \cite{fukunaga2022block}, the authors consider a semi-unbalanced problem in which the penalty is imposed only on one of the two marginals \cite{liu2026solving}.
To the best of our knowledge, most of the works within this field are based on approaching the problem from the dual point of view and by adding a regularization term to the objective function, be it the Kullback Leibler divergence \cite{chizat2018scaling} or be it an $\ell_2$ penalization term \cite{nguyen2023unbalanced}.
A relatively under-explored research direction is to leverage some well-established algorithm for constrained optimization, such as the Frank-Wolfe (FW) algorithm, to retrieve the optimal solution for the primal formulation.
Indeed, to the best of our knowledge, the only papers that use the FW routine to retrieve the solution to the UOT problem are \cite{sejourne2022faster}, where the authors use the FW routine to approach the dual formulation induced by the one dimensional UOT problem, and \cite{fukunaga2022block}, where the authors define a block-coordinate version of the algorithm that exploits the structure of their semi-unbalanced problem.
Both papers showed that it is possible to define fast and efficient Frank-Wolfe routines for the UOT whenever one can harness the structure of the problem: in \cite{sejourne2022faster} they do so by showing that the Linear Minimization Oracle associated with each iteration can be performed in closed form, while in \cite{fukunaga2022block} it is done by decomposing the feasible set into the direct product of small simplex sets.
{\color{black}
The Frank-Wolfe algorithm, also known as the conditional gradient method, is a classic optimization algorithm for constrained problems introduced by Frank and Wolfe~\cite{frank1956algorithm}.
Classical convergence guarantees have been extensively studied \cite{canon1968tight} along with the conditions that guarantee linear convergence \cite{levitin1966constrained,demianovapproximate} and enable to identify the solution support in a finite number of steps \cite{bomze2020active}.
To ease the conditions needed to make the standard FW algorithm converge linearly, several variants, including away-step \cite{wolfe1970convergence}, pairwise \cite{demianovapproximate}, and fully corrective Frank-Wolfe, have been proposed to improve convergence rates under suitable geometric assumptions~\cite{lacoste2015global}.
In recent years, the interest for the Frank-Wolfe algorithm and its variants has renewed due to its applications in sub-modular optimization problems \cite{bach2013learning}, variational inference problems \cite{krishnan2015barrier}, traffic routing \cite{fukushima1984modified}, clustering \cite{bomze1997evolution}, and sparse neural network training \cite{grigas2019stochastic}.
We refer the reader to \cite{bomze2021frank} and \cite{braun2022conditional} for a survey on the family of Frank-Wolfe algorithms.
}

\subsection{Our Contribution}
In this paper, we propose a novel approach to solve the Unbalanced Optimal Transportation (UOT) problem presented in \cite{liero2018optimal} that hinges on a clever exploitation of the structure of the UOT solutions.
Our paper is structured as follows.
First, we study the properties of the solution to the UOT problem.
In particular, given any couple of discrete positive measures $\mu$ and $\nu$, we show that there exists at least an optimal transportation plan between $\mu$ and $\nu$ that is discrete and characterize the finite set of points on which such solution is supported.
Building on this preliminary result, we show that any solution has finite measure and provide an upper and lower bound on the marginals of any solution that depend only on the specific of the problem.
In particular, we provide an upper bound on the measure that any solution assigns to any couple of points in its support.
These upper bounds reveal that, under suitable assumptions, the maximum distance between couples of points that belong to the support is bounded by a global constant depending on the cost function and the Csiszár divergence defining the problem.
While a similar characterization was known for the Total Variation case \cite{sejourne2023unbalanced}, our approach can be applied to any UOT problem following the formulation presented in \cite{liero2018optimal}.
We then propose a truncated version of the UOT problem that is obtained by introducing a parameter $R$ that saturates the cost function $c$.
We show that the truncated problem admits a simple formulation, reducing the total number of variables needed to describe the associated UOT problem.
Moreover, any truncated problem naturally yields a feasible point for the original UOT problem that acts as an approximated solution of the problem.
We provide an upper bound on the primal sub-optimality gap between the real solution and the approximated solution retrieved from the truncated formulation, revealing that the gap decreases as a polynomial function of $R$.
Noticeably, all of our proposed bounds on the optimal solution apply to general framework.

We then specialize our study to the case in which the Csiszár divergence belongs to the class of powerlike entropies, a class of entropies that includes the Kullback Leibler and the quadratic penalty function as specific sub-cases.
To approach this class of problems, we define an efficient implementation of the Pairwise Frank-Wolfe routine capable of solving the UOT problem for any powerlike entropy and any couple of positive measures.
Our routine cleverly combines our theoretical results to save the iterates of the problem in a fast and efficient way.
As a by product, we also study the convergence speed of the Pairwise Frank-Wolfe algorithm under the assumptions that the objective has local Hölder continuous gradient and is $p$-uniformly convex.
Moreover, we provide a support identifiability property for the Pairwise Frank-Wolfe algorithm, which ensures that our routine finds the support of the optimal solution within a finite number of steps.

\section{Notation and problem formulation}
\label{sec:notation}

We collect the notation and assumptions used throughout the paper.

\subsection{Basic notation}

We write $\R_+=\left[0,+\infty\right[$, $\R_{++}=\left]0,+\infty\right[$, and $[d]=\{1,\ldots,d\}$. 
The class $\Gamma_0(\R^d)$ consists of proper, lower semicontinuous convex functions from $\R^d$ to $\left]-\infty,+\infty\right]$.
For a set $C$, $\ind_C$ denotes its indicator function, equal to zero on $C$ and $+\infty$ elsewhere. 
We use $\partial f$ for the convex subdifferential and $N_C=\partial\ind_C$ for the normal cone of a closed convex set. 
Interior and relative interior are denoted by $\operatorname{int}$ and $\ri$, respectively.

For a finite-dimensional vector $z$, $\|z\|_p$ is its usual $\ell_p$ norm; $\|z\|$ denotes the Euclidean norm and $\|z\|_\infty:=\max_{i\in[d]}|z_i|$ denotes the maximum norm. 
Under these notations, we set $B_{h,p}(z)=\{w:\|w-z\|_p\leq h\}$.
The inner product on matrices is $\pair{a}{b}=\sum_{i,j}a_{ij}b_{ij}$, and $\one_d$ is the vector of $d$ ones.
For a nonnegative vector or matrix, $\|z\|_1$ is the sum of its entries.
For a real number $t$, put $t_+=\max\{t,0\}$.
For a compact metric space $\ZC$, $\Mplus(\ZC)$ is the set of finite nonnegative Borel measures and $\CC(\ZC)$ is the space of continuous real functions. 
Given $\ZC^\prime\subseteq\ZC$ and $\mu\in\Mplus(\ZC)$, we denote by $\mu_{|\ZC^\prime}$ the restriction of $\mu$ to $\ZC^\prime$, that is $\mu_{|\ZC^\prime}(A)=\mu(A\cap\ZC^\prime)$ for every measurable $A\subseteq\ZC$.
The support of a measure $\mu$ is denoted by $\supp(\mu)$.
In what follows, we say that a measure is a discrete measure if it has finite support. 
If $T$ is measurable, $T_\#\mu(E)=\mu(T^{-1}(E))$ denotes the pushforward. 
For a measure $\pi\in\Mplus(\XC\times\YC)$, its marginals are $\pi^{(1)}=(\operatorname{pr}_\XC)_\#\pi$ and $\pi^{(2)}=(\operatorname{pr}_\YC)_\#\pi$.
Let $f$ be real-valued and differentiable on an open set containing $\Omega\subset\R^d$.
Its gradient is $\alpha$-H\"older continuous on $\Omega$, with $\alpha\in\left]0,1\right]$, if there is a finite $H_\alpha\geq0$ such that
\begin{equation*}
 \|\nabla f(x)-\nabla f(y)\|\leq H_\alpha\|x-y\|^\alpha
 \qquad(x,y\in\Omega).
\end{equation*}
The case $\alpha=1$ is Lipschitz continuity of the gradient.
For $p\geq2$, $f$ is $p$-uniformly convex on $\Omega$ with constant $\kappa_p>0$ if
\begin{equation}
 f(y)\geq f(x)+\pair{\nabla f(x)}{y-x}
               +\frac{\kappa_p}{p}\|y-x\|^p
 \qquad(x,y\in\Omega).
 \label{eq:growth_cond_intro}
\end{equation}
For $p=2$, this is strong convexity. 
Notice that these constants may depend on
$\Omega$; a global assertion requires the estimate on the whole set under consideration with one constant.
Finally, if $f$ has nonempty minimizer set $\XC^\star$ and minimum $f^\star$, it is $p$-conditioned on $\Omega\supset \XC^\star$ if there is a constant $\upsilon_p>0$ such that
\begin{equation}
\label{eq:pcond_intro}
 f(x)-f^\star\geq\upsilon_p\operatorname{dist}(x,\XC^\star)^p
 \qquad(x\in\Omega),
\end{equation}
where $\operatorname{dist}(x,\XC^\star):=\min_{x^\star\in\XC^\star}\|x-x^\star\|$ is the standard element-to-set distance.

\subsection{Standing assumptions, marginal penalties, and problem formulation}
\label{subsec:setting}

Let $\XC$ and $\YC$ be nonempty compact metric spaces, let $c : \XC\times \YC\to\left[0,+\infty\right[$ be continuous, and let $\mu\in\Mplus(\XC)$ and $\nu\in\Mplus(\YC)$ be finite nonnegative measures.
Given $x\in\XC$, we denote by $\delta_x\in\Mplus(\XC)$ the Dirac delta centred in $x$.
In what follows, we always assume that the support of $\mu$ and $\nu$ are finite so that $\mu$ and $\nu$ are positive linear combination of different centred Dirac deltas.
Moreover, $\phi\colon\R\to\left]-\infty,+\infty\right]$ is proper, lower semicontinuous, and convex, and satisfies
\begin{equation}
 \phi\geq0,\qquad ]0,+\infty[ \subseteq\dom(\phi)\subseteq[0,+\infty[,
 \qquad \text{and} \qquad 
 0\in\partial\phi(1),
 \label{eq:assumptions}
\end{equation}
where $\dom(\phi)$ denotes the domain of $\phi$.
In particular, $\phi$ is finite on $]0,+\infty[$, whereas $\phi(0)$ may be infinite. 
In what follows, we assume that
\begin{equation}
 \phipi:=\lim_{t\to+\infty}\frac{\phi(t)}{t}\in\left]0,+\infty\right],
 \qquad 
 \phipz:=\lim_{t\downarrow0}\phi^\prime_+(t)\in\left[-\infty,0\right].
 \label{eq:slopes}
\end{equation}
When $\phi(0)<+\infty$, $\phipz$ is the right derivative of $\phi$ at zero; when $\phi(0)=+\infty$, $\phipz=-\infty$.
For $t>0$, let $\phi^\prime_+(t)$ and $\phi^\prime_-(t)$ denote the right and left derivatives at $t$, then
\begin{equation}
 \partial\phi(t)=[\phi'_-(t),\phi'_+(t)],
 \qquad 
 \phipz\leq\phi'_-(t)\leq\phi'_+(t)\leq\phipi.
 \label{eq:derivative-order}
\end{equation}
At zero, $\partial\phi(0)=]-\infty,\phipz]$ if $\phipz$ is finite, and $\partial\phi(0)=\varnothing$ otherwise.
For a reference $m\geq0$, define the closed marginal penalty
\begin{equation}
 \varphi_m(t):=
 \begin{cases}
  m\phi(t/m),&m>0,\ t\geq0,\\
  \phipi t,&m=0,\ t\geq0,\ \phipi<+\infty,\\
  0,&m=0,\ t=0,\ \phipi=+\infty,\\
  +\infty,&\text{otherwise}.
 \end{cases}
 \label{eq:marginal-penalty}
\end{equation}
Thus $\varphi_0=\ind_{\{0\}}$ when $\phipi=+\infty$. 
The function in \eqref{eq:marginal-penalty} is also known as the perspective of $\phi$ with the reference $m$ fixed. 
By definition, the subdifferential of $\varphi_m$ is
\begin{equation}
 \partial \varphi_m(t)=
 \begin{cases}
  \partial\phi(t/m),&m>0,\\
  \{\phipi\},&m=0,\ \phipi<+\infty,\ t>0,\\
  ]-\infty,\phipi],&m=0,\ \phipi<+\infty,\ t=0,\\
  \R,&m=0,\ \phipi=+\infty,\ t=0,\\
  \varnothing,&\text{otherwise}.
 \end{cases}
 \label{eq:marginal-subdifferential}
\end{equation}
Under these assumptions, the following increment bound holds
\begin{equation}
 \varphi_m(x+\delta)-\varphi_m(x)\leq\phipi\delta \qquad(\phipi<+\infty,\ x\in\dom (\varphi_m),\ \delta\geq0).
 \label{eq:upper-increment}
\end{equation}
For $\delta=0$ the bound is trivial.
If $m>0$ and $\delta>0$, it follows from the convexity of $\phi$
\begin{equation}
\label{eq:perspective_bound}
 \frac{\varphi_m(x+\delta)-\varphi_m(x)}{\delta}  =\frac{\phi((x+\delta)/m)-\phi(x/m)}{\delta/m}
 \leq\phi'_-((x+\delta)/m)\leq\phipi,
\end{equation}
while for $m=0$ the bound holds as an equality since $\varphi_0$ is a linear function if $\phipi<\infty$.

Given $\alpha,\mu\in\Mplus(\XC)$ two positive measures, let $\alpha=\rho\mu+\alpha^\perp$ be the unique Lebesgue decomposition of $\alpha$ with respect to $\mu$.
The Csiszár divergence induced by $\phi$ is
\begin{equation}
\label{def:Dphi}
 \Dphi(\alpha\mid\mu) =\int_\XC\phi(\rho)\,d\mu+\phipi\alpha^\perp(\XC),
\end{equation}
where we adopted the convention $+\infty\cdot0=0$. 
Then, given $\mu\in\Mplus(\XC)$ and $\nu\in\Mplus(\YC)$, the Unbalanced Optimal Transport (UOT) problem is defined as
\begin{equation}
 v_c:=\inf_{\pi\in\Mplus(\XC\times \YC)}\UTphi(\pi),\qquad
 \UTphi(\pi):=\int_{\XC\times \YC}c\,d\pi
       +\Dphi(\pi^{(1)}\mid\mu)
       +\Dphi(\pi^{(2)}\mid\nu).
 \label{eq:UOTphi}
\end{equation}

\begin{remark}[Zero input measures]
\label{rem:zero-measures}
If both input measures are null measures, either the zero measure is the unique finite-valued minimizer of $\UTphi$ or $\UTphi(\pi)=+\infty$ for every $\pi\in\Mplus(\XC\times\YC)$.
If $\phipi=+\infty$ and exactly one input measure is the null measure, every nonzero plan has infinite objective value.
In this case, the zero plan then has value $\phi(0)(\mu(\XC)+\nu(\YC))$.
In particular, it is the unique finite-valued minimizer when $\phi(0)<+\infty$ while when $\phi(0)=+\infty$ there are no finite-valued plans.
For this reason, we do not consider these cases and hence we assume
\begin{equation*}
 \begin{cases}
 \mu(\XC)+\nu(\YC)>0,&\phipi<+\infty,\\
 \mu(\XC)>0\; \text{and}\; \nu(\YC)>0,&\phipi=+\infty.
 \end{cases}
\end{equation*}
\end{remark}

We recall that, under the stated assumptions, the Unbalanced Optimal Transport problem always admits at least a solution \cite{liero2018optimal}.

\begin{proposition}[Existence]
\label{prop:existence}
Under the nondegenerate conditions in Remark~\ref{rem:zero-measures}, $v_c<+\infty$ and the infimum in \eqref{eq:UOTphi} is attained.
\end{proposition}

\section{The UOT Problem between Discrete Positive Measures}

We first show that the UOT problem induced by any couple of positive measures with finite support admits a minimizer with finite support.
This allows to cast the UOT problem to a finite-dimensional convex problem and to express the optimality conditions via potential vectors.
We then exploit the optimality conditions to derive bounds on the support, the individual entries, and the marginals of optimal solutions to the UOT problem.

\subsection{The Support of the Optimal Solution}
\label{subsec:finite-support}

In what follows, we set $\AC=\supp(\mu)$ and $\BC=\supp(\nu)$.
If $\phipi=+\infty$, finiteness of the divergences forces every plan attaining bounded objective value to be supported on $\AC\times \BC$ \cite{liero2018optimal}.
In the next result, we show that something similar happens when $\phipi<+\infty$.
Indeed, let us set
\begin{equation}
 \hat y(x)\in\argmin_{y\in \YC}c(x,y)\quad(x\in \AC),\qquad
 \hat x(y)\in\argmin_{x\in \XC}c(x,y)\quad(y\in \BC),
 \label{eq:selected-partners}
\end{equation}
and define
\begin{equation}
 \KC_\XC=\AC\cup\{\hat x(y):y\in \BC\},\qquad
 \KC_\YC=\BC\cup\{\hat y(x):x\in \AC\}.
 \label{eq:enlarged-grids}
\end{equation}
The selected points exist by compactness and continuity, albeit they might not be unique, in which case we pick an arbitrary minimizer.

\begin{theorem}
\label{thm:solution_support_characterization}
    Let $c:\XC\times \YC\to [0,+\infty[$ be a continuous cost function and $\phi\in\Gamma_0(\R)$ be such that $0<\phipi<+\infty$.
    Then, for any couple of finite positive measures $\mu\in\Mplus(\XC)$ and $\nu\in\Mplus(\YC)$, every solution to the UOT, namely $\pi$, satisfies
    \begin{equation}
     \pi(\AC^c\times \BC^c)=0 \quad \quad \text{with}\;\; \AC=\supp(\mu)\quad\text{and}\quad \BC=\supp(\nu).
     \label{eq:no-double-singular}
    \end{equation}
    Moreover, there exists at least an optimal solution to the UOT problem, namely $\pi^\prime$, such that
    \begin{equation}
        \label{eq:set_spt_opt_plan}
        \supp(\pi^\prime)\subseteq\KC_\XC\times\KC_\YC,
    \end{equation}
    where $\KC_\XC$ and $\KC_\YC$ are defined as in \eqref{eq:enlarged-grids}.
\end{theorem}

\begin{proof}
    First, we show that $\pi(\AC^c\times \BC^c)=0$.
    Indeed, towards a contradiction, assume that $\pi$ is a solution to \eqref{eq:UOTphi} such that $\pi\big(\AC^c\times \BC^c\big)=\lambda>0$.
    We then define $\pi^\prime=\pi-\pi_{|\AC^c\times \BC^c}$.
    Notice that $\pi^\prime$ is a positive measure, thus feasible for \eqref{eq:UOTphi}.
    We now show that $\UTphi(\pi^\prime)<\UTphi(\pi)$.
    First, since $\pi_{|\AC^c\times \BC^c}$ is a positive measure and $c$ is non-negative 
    \begin{equation*}
        \int_{\XC\times\YC} c(x,y) \pi^\prime(dx, dy) \le \int_{\XC\times\YC} c(x,y) \pi(dx, dy).
    \end{equation*}
    Moreover, the first and second marginals of $\bar\pi:=\pi_{|\AC^c\times \BC^c}$ are, by construction, singular with respect to $\mu$ and $\nu$, respectively.
    Thus
    \begin{align*}
        \pi^{(1)}:=\rho^{(\mu)}\mu+\alpha_\perp\quad &\Rightarrow\quad (\pi^\prime)^{(1)}=\pi^{(1)}-\bar\pi^{(1)}=\rho^{(\mu)}\mu+\alpha_\perp-\bar\pi^{(1)}\\
        \pi^{(2)}:=\rho^{(\nu)}\nu+\beta_\perp \quad &\Rightarrow\quad (\pi^\prime)^{(2)}=\pi^{(2)}-\bar\pi^{(2)}=\rho^{(\nu)}\nu+\beta_\perp-\bar\pi^{(2)},
    \end{align*}
    and, by the uniqueness of Lebesgue's decomposition for measures, we infer $\alpha_\perp^\prime :=\alpha_\perp-\bar\pi^{(1)}$ and $\beta_\perp^\prime :=\beta_\perp-\bar\pi^{(2)}$ are the singular part of the first and second marginals of $\pi^\prime$ with respect to $\mu$ and $\nu$, respectively.
    Therefore, the definition of $\Dphi$ yields
    \begin{equation}
        \Dphi\big((\pi^{\prime})^{(1)}|\mu\big)+\Dphi\big((\pi^{\prime})^{(2)}|\nu\big) = \Dphi(\pi^{(1)}|\mu)+\Dphi(\pi^{(2)}|\nu) -2\lambda\phipi,
    \end{equation}
    and hence $\UTphi(\pi^\prime)<\UTphi(\pi)$, which is a contradiction allowing us to conclude \eqref{eq:no-double-singular}.
    We now move to the second part of the statement.
    Without loss of generality, we assume that $\pi(\AC^c\times \BC^c)=0$, which implies
    \begin{equation}
    \label{eq:suppsupsetc}
        \supp(\pi)^c\supseteq \AC^c\times \BC^c.
    \end{equation}
    Given $\pi$ a solution to \eqref{eq:UOTphi}, owing \eqref{eq:suppsupsetc} we have
    \begin{equation}
    \label{eq:supppic_part}
        \supp(\pi)\backslash (\KC_\XC\times \KC_\YC) \subseteq \big(\mathcal{S}_0\times\KC_\YC^c\big)
        \cup\big(\KC_\XC^c\times\mathcal{T}_0\big), \ \text{with}\  \mathcal{S}_0=\AC,\;\mathcal{T}_0=\BC.
    \end{equation}
    We note that the two sets on the right hand side of \eqref{eq:supppic_part} are disjoints, since $\KC_\XC^c$ is disjoint from $\AC$ by definition and, similarly, $\KC_\YC^c$ is disjoint from $\BC$.
    If the left hand side of \eqref{eq:supppic_part} is empty, we have nothing to prove since $\supp(\pi)\backslash(\KC_\XC\times\KC_\YC)=\emptyset$ implies $\supp(\pi)\subseteq \KC_\XC\times\KC_\YC$.
    Suppose therefore that there exists an element $(\bar x,\bar y)\in\supp(\pi)\backslash (\KC_\XC\times \KC_\YC)$ and, without loss of generality, that $(\bar x,\bar y)\in\mathcal{S}_0\times \KC_\YC^c$.
    We will construct a transport plan $\tilde \pi$ which is still optimal for \eqref{eq:UOTphi} and satisfies
    \begin{equation}
    \label{eq:11092026}
        \tilde \pi_{\lvert (\{\bar x\}\times \YC)^c} = \pi_{\lvert (\{\bar x\}\times \YC)^c} \quad\text{and}\quad \supp(\tilde \pi_{\lvert \{\bar x\}\times \YC})\subseteq \KC_{\XC}\times\KC_{\YC},
    \end{equation}
    so that
    \begin{equation}
    \label{eq:supppic_part2}
        \supp(\tilde\pi)\backslash (\KC_\XC\times \KC_\YC) \subseteq \big(\mathcal{S}_1\times\KC_\YC^c\big)\cup\big(\KC_\XC^c\times\mathcal{T}_0\big), \ \text{with}\  \mathcal{S}_1=\mathcal{S}_0\backslash\{\bar x\}.
    \end{equation}
    Similarly, if $(\bar x,\bar y)\in\KC^c_{\XC}\times \mathcal{T}_0$, then there exists an optimal transport plan $\tilde \pi$ such that $\supp(\tilde\pi)\backslash (\KC_\XC\times \KC_\YC) \subseteq \big(\mathcal{S}_0\times\KC_\YC^c\big)\cup\big(\KC_\XC^c\times\mathcal{T}_1\big)$ with $\mathcal{T}_1=\mathcal{T}_0\backslash\{\bar y\}$.
    If starting from a solution whose support satisfies \eqref{eq:supppic_part} we are able to build a solution that satisfies \eqref{eq:supppic_part2}, we can iteratively apply this construction.
    Since at each iteration we reduce the cardinality of either $\mathcal{S}_0$ or $\mathcal{T}_0$ by one, within a finite number of iteration, we obtain a solution to \eqref{eq:UOTphi}, namely $\hat \pi$, such that $\supp(\hat \pi)\subseteq \KC_\XC\times\KC_\YC$, proving the statement.
    In order to complete the proof, we have to prove that if $\pi$ is optimal, satisfies
    \begin{equation}
    \label{eq:supppic_part3}
        \supp(\pi)\backslash (\KC_\XC\times \KC_\YC) \subseteq \big(\mathcal{S}\times\KC_\YC^c\big)
        \cup\big(\KC_\XC^c\times\mathcal{T}\big), \ \text{with}\  \mathcal{S}\subseteq\supp(\mu),\;\mathcal{T}\subseteq\supp(\nu),
    \end{equation}
    and $(\bar x,\bar y)\in\mathcal{S}\times \KC_\YC^c$, then there exist an optimal $\tilde \pi$ that satisfies \eqref{eq:11092026}.
    By disintegrating $\pi$ with respect to $\pi^{(1)}$, we find a measurable family of measures, namely $\{\rho(\cdot|x)\}_{x\in\XC}$, such that for every pair of measurable sets $X\subset\XC$ and $Y\subset\YC$, we get
    \begin{align*}
        \pi(X\times Y) &=\int_{X}\rho(Y|x)\pi^{(1)}(dx)=\int_{X}\rho(Y|x)f(x)\mu(dx)+\int_{X}\rho(Y|x)\pi^{(1)}_{\perp,\mu}(dx)
    \end{align*}
    where in the last equality we used the standard Lebesgue's decomposition of the positive measures $\pi^{(1)}=f\mu+\pi^{(1)}_{\perp,\mu}$.
    For every $Y\subseteq\YC$, we define
    \begin{equation}
        \tilde\rho(Y|x)=\begin{cases}
            \rho(\BC^c|\bar x)\delta_{\hat y(\bar x)}(Y)+\rho(\BC\cap Y|\bar x) \quad\quad &\text{if}\;\; x=\bar x;\\
            \rho(Y|x) \quad\quad &\text{otherwise},
        \end{cases}
    \end{equation}
    where $\hat{y}(\bar x)\in\argmin_{y\in \YC}c(\bar x,y)$ and set $\tilde\pi=\tilde\rho\otimes \pi^{(1)}$.
    First, we show that $\tilde\pi^{(1)}=\pi^{(1)}$.
    Indeed, given a measurable set $X\subseteq \XC$, we have
    \begin{align*}
        \tilde \pi^{(1)}(X) &= \tilde \pi\big(X\times \YC\big) =\int_{X}\tilde\rho\big(\YC|x) \pi^{(1)}(dx) = \pi^{(1)}(X),
    \end{align*}
    since $\tilde\rho(\YC|x)=\rho(\YC|x)=1$ for every $x\in \XC$.
    Notice that $\tilde\pi^{(1)}=\pi^{(1)}$ implies $\Dphi(\tilde\pi^{(1)}|\mu)=\Dphi(\pi^{(1)}|\mu)$ and $\tilde\pi(\XC\times \YC)=\pi(\XC\times \YC)$.
    By construction, $\tilde \pi$ satisfies \eqref{eq:11092026}.
    Indeed, given $C\subseteq \{\bar x\}^c\times \YC$, let $\chi_C:\XC\times\YC\to\R$ denote the function that takes value $1$ when $(x,y)\in C$ and $0$ otherwise, then
    \begin{align*}
        \tilde\pi(C)&=\int_{\XC\times\YC}\chi_C(x,y)\tilde\pi(dx,dy)=\int_{\XC\backslash\{\bar x\}}\int_{\YC}\chi_C(x,y)\tilde\rho(dy|x)\tilde\pi^{(1)}(dx)\\
        &=\int_{\XC\backslash\{\bar x\}}\int_{\YC}\chi_C(x,y)\rho(dy|x)\pi^{(1)}(dx)=\int_{\XC}\int_{\YC}\chi_C(x,y)\rho(dy|x)\pi^{(1)}(dx)=\pi(C),
    \end{align*}
    since $C\subseteq \{\bar x\}^c\times \YC$ and $\tilde\rho(\cdot|x)=\rho(\cdot|x)$ for every $x\neq \bar x$.
    It then follows that $\tilde\pi_{|(\{\bar x\}\times \YC)^c}=\pi_{|(\{\bar x\}\times \YC)^c}$.
    Moreover, $\supp(\tilde\pi_{|\{\bar x\}\times\YC})\subset\KC_\XC\times\KC_\YC$ since $\KC_\XC\times\KC_\YC$ is closed and
    \begin{align*}
       \tilde\pi_{|\{\bar x\}\times\YC}\big((\KC_\XC\times\KC_\YC)^c\big) &\le\tilde\pi\big((\{\bar x\}\times\YC)\cap(\KC_\XC^c\times\mathcal{Y})\big)+\tilde\pi\big((\{\bar x\}\times\YC)\cap(\XC\times\KC_\YC^c)\big)\\
        &= \tilde\pi\big(\{\bar x\}\times\KC_\YC^c\big) =\int_{\{\bar x\}}\tilde \rho(\KC_\YC^c|\bar x)\tilde \pi^{(1)}(dx)=0,
    \end{align*}
    where in the second last equality we used the fact that $\KC_\YC^c$ is disjoint from $\BC\cup\{\hat y(x)\,:\, x\in \AC\}$ and in the last equality we used the definition of $\tilde \rho$.
    We now consider $\tilde\pi^{(2)}$.
    Given a measurable $Y\subseteq \YC$, we have
    \begin{align}
        \nonumber\tilde \pi^{(2)} (Y) & 
        = \int_{\XC}\tilde\rho(Y|x)\pi^{(1)}(dx) = \int_{\XC}\tilde\rho(Y|x)f(x)\mu(dx)+\int_{\XC}\tilde\rho(Y|x)\pi^{(1)}_{\perp,\mu}(dx)\\
        \label{eq:form_pi_tild_2}&=f_{\bar x}\mu_{\bar x}\big(\rho(\BC^c|\bar x)\delta_{\hat y(\bar x)}(Y) + \rho(Y \cap \BC|\bar x)\big)\\
        \nonumber&\quad +\sum_{x\neq \bar x,\ x\in \AC}\rho(Y|x)f_x\mu_x +\int_{\XC}\rho(Y|x)\pi^{(1)}_{\perp,\mu}(dx),
    \end{align}
    where we have used the definition of $\tilde \rho$.
    Since $\nu$ is discrete
    \begin{equation}
        \tilde\pi^{(2)}=\tilde\pi^{(2)}_\nu+\tilde\pi^{(2)}_{\perp,\nu}, \quad\quad \text{with}\quad \tilde\pi^{(2)}_\nu=\sum_{y\in\BC}\frac{\tilde\pi^{(2)}(\{y\})}{\nu_y}\nu_y\delta_{y},\;\; \nu_y=\nu(\{y\}),
    \end{equation}
    and $\tilde\pi^{(2)}_{\perp,\nu}(\{y\})=0$ for every $y\in\BC$.
    To conclude, we show $\UTphi(\tilde\pi)\le\UTphi(\pi)$ by proving that
    \begin{equation}
    \label{eq:upperboundtranspcost_thm_spt}
        \int_{\XC\times \YC}c(x,y)\tilde\pi(dx,dy)\le\int_{\XC\times \YC}c(x,y)\pi(dx,dy)
    \end{equation}
    and $\Dphi(\tilde\pi^{(2)}|\nu)\le \Dphi(\pi^{(2)}|\nu)$ since we have already proven that $\Dphi(\tilde\pi^{(1)}|\mu)\le \Dphi(\pi^{(1)}|\mu)$.
    Let us start with \eqref{eq:upperboundtranspcost_thm_spt}.
    By definition,
    \begin{align*}
        \int_{\XC\times \YC}c(x,y)\pi(dx,dy)&=\int_{\XC\backslash\{\bar x\}}\bigg(\int_\YC c(x,y)\rho(dy|x)\bigg)\pi^{(1)}(dx)+\pi^{(1)}(\{\bar x\}) \int_\YC c(\bar x,y)\rho(dy|\bar x).
    \end{align*}
    Notice that
    \begin{align*}
        \int_\YC c(\bar x,y)\rho(dy|\bar x)\ge c(\bar x,\hat y(\bar x))\rho(\BC^c|\bar x)+\sum_{y\in \BC}c(\bar x,y)\rho(\{y\}|\bar x) = \int_\YC c(\bar x,y)\tilde\rho(dy|\bar x).
    \end{align*}
    To conclude, we recall that $\tilde\rho(dy|x)=\rho(dy|x)$ for every $x\neq \bar x$, hence \eqref{eq:upperboundtranspcost_thm_spt} follows.
    Let us now consider $\Dphi(\tilde\pi^{(2)}|\nu)\le \Dphi(\pi^{(2)}|\nu)$.
    We split our discussion in two parts, depending on whether $\hat y(\bar x)\in \BC$ or $\hat y(\bar x)\notin \BC$.
    Let us first consider the case in which $\hat y(\bar x)\notin \BC$.
    In this case, for every $y\in \BC$ we have
    \begin{align}
    \label{eq:equality_snu}
        \tilde\pi^{(2)}(\{y\})&=\sum_{x\in \AC}\rho(\{y\}|x)f_x\mu_x +\int_{\XC}\rho(\{y\}|x)\pi^{(1)}_{\perp,\mu}(dx)=\pi^{(2)}(\{y\}),
    \end{align}
    hence $\tilde\pi^{(2)}(\{y\})=\pi^{(2)}(\{y\})$ and thus $\tilde\pi^{(2)}$ and $\pi^{(2)}$ have the same absolute continuous part with respect to $\nu$.
    Concerning the singular parts of $\tilde\pi^{(2)}$ and $\pi^{(2)}$, we note that by \eqref{eq:form_pi_tild_2}, the following holds 
    \begin{align*}
        \tilde\pi^{(2)}(\BC^c)&=\rho(\BC^c|\bar x)f_{\bar x}\mu_{\bar x}+\sum_{x\neq \bar x,\ x\in \AC}\rho(\BC^c|x)f_x\mu_x +\int_{\XC}\rho(\BC^c|x)\pi^{(1)}_{\perp,\mu}(dx)\\
        & = \sum_{x\in \AC}\rho(\BC^c|x)f_x\mu_x +\int_{\XC}\rho(\BC^c|x)\pi^{(1)}_{\perp,\mu}(dx) = \pi^{(2)}(\BC^c).
    \end{align*}
    Therefore $\tilde\pi^{(2)}_{\perp,\nu}(\YC)=\pi^{(2)}_{\perp,\nu}(\YC)$, and hence, in the end, $\Dphi(\tilde\pi^{(2)}|\nu)=\Dphi(\pi^{(2)}|\nu)$.
    Let us now consider the case in which $\hat y (\bar x)\in \BC$, so that $\hat y(\bar x)=\bar y$ for some $\bar y\in \BC$. 
    For every $y\neq \hat y(\bar x)$ equation \eqref{eq:equality_snu} still holds. 
    If $y=\bar y$
    \begin{align*}
        \tilde\pi^{(2)}(\{\bar y\})&= \big(\rho(\BC^c|\bar x)+\rho(\{\bar y\}|\bar x)\big)f_{\bar x}\mu_{\bar x} + \sum_{x\neq\bar x}\rho(\{\bar y\}|x)f_x\mu_x +\int_{\XC}\rho(\{\bar y\}|x)\pi^{(1)}_{\perp,\mu}(dx)\\
        &= \pi^{(2)}(\{\bar y\}) + \rho(\BC^c|\bar x)f_{\bar x}\mu_{\bar x}.
    \end{align*}
    Let us set $\bar\Delta=\rho(\BC^c|\bar x)f_{\bar x}\mu_{\bar x}$.
    Since $\tilde\pi(\XC\times \YC)=\pi(\XC\times \YC)$, we have
    \begin{equation}
        \tilde \pi^{(2)}(\BC^c)=\pi^{(2)}(\BC^c)- \bar\Delta.
    \end{equation}
    Lastly, if $\bar \Delta=0$, we have $\tilde\pi^{(2)}=\pi^{(2)}$ and thus $\Dphi(\tilde\pi^{(2)}|\nu)=\Dphi(\pi^{(2)}|\nu)$.
    If $\bar\Delta\neq0$, from \eqref{eq:perspective_bound} we infer
    \begin{align*}
        \Dphi(\tilde\pi^{(2)}|\nu)-\Dphi(\pi^{(2)}|\nu)&=\nu_{\bar y}\bigg(\phi\Big(\frac{\pi^{(2)}_{\bar y}+\bar\Delta}{\nu_{\bar y}}\Big)-\phi\Big(\frac{\pi^{(2)}_{\bar y}}{\nu_{\bar y}}\Big)\bigg)-\phipi\bar\Delta\le 0.
    \end{align*}
    In both cases, we have $\Dphi(\tilde\pi^{(2)}|\nu)\le\Dphi(\pi^{(2)}|\nu)$, which concludes the proof.
\end{proof}

\begin{remark}[Existence versus all optimizers]
\label{rmk:supp_prod}
If $\phipi<+\infty$, Theorem \ref{thm:solution_support_characterization} ensures the existence of a discrete optimal solution, not discreteness of every optimal solution.
If all the minima in \eqref{eq:selected-partners} are attained at a unique point, the proof of Theorem \ref{thm:solution_support_characterization} shows that every optimal plan is supported on $\KC_\XC\times \KC_\YC$.
In particular, if $\XC=\YC$ and $c(x,y)=0$ if and only if $x=y$, then $\KC_\XC=\KC_\YC=\AC\cup \BC$ and every optimal solution is supported on $(\AC\cup \BC)\times(\AC\cup \BC)$.
\end{remark}

Owing to Theorem \ref{thm:solution_support_characterization}, given two discrete positive measures $\mu\in\Mplus(\XC)$ and $\nu\in\Mplus(\YC)$, there exists a finite set of points $\{x_1, \dots, x_n\} \subseteq \XC$ and $\{y_1, \dots, y_m\}\subseteq \YC$ such that there exist a solution to \eqref{eq:UOTphi} for which it holds
\begin{equation*}
    \supp(\pi)\subseteq \{x_1, \dots, x_n\}\times \{y_1, \dots, y_m\}.
\end{equation*}
In particular, it is possible to formulate the UOT problem as a minimization problem over the set of $n\times m$ matrices with positive entries, that is
\begin{equation}
 v_c=\min_{\pi\in\R_+^{n\times m}}F_c(\pi),\qquad
 F_c(\pi)=\pair{c}{\pi}
        +\sum_{i=1}^n \varphi_{\mu_i}\big((\Pruno\pi)_i\big)
        +\sum_{j=1}^m \varphi_{\nu_j}\big((\Prdue\pi)_j\big),
 \label{eq:matrix-problem}
\end{equation}
where $\Pruno\colon \R^{n\times m}\to \R^n$ and $\Prdue\colon \R^{n\times m}\to \R^m$ are the linear maps corresponding to the operations of push-forward of the projections onto $\XC$ and $\YC$, respectively, that is
\begin{equation*}
    \Pruno \pi = \bigg(\sum_{j=1}^m \pi_{ij}\bigg)_{1 \leq i\leq n}
    \quad\text{and}\quad
    \Prdue \pi = \bigg(\sum_{i=1}^n \pi_{ij}\bigg)_{1 \leq j\leq m},
\end{equation*}
and $\varphi_m$ is the perspective function defined as in \eqref{eq:marginal-penalty}.
Even if the problem can be formulated via matrices, we will refer to any $\pi$ feasible for \eqref{eq:matrix-problem} as a transportation plan.
With a slight abuse of notation, we will refer to $\pi_{ij}$ as the entry of the transportation plan.
Notice that in \eqref{eq:matrix-problem} all reference weights $\{\mu_i\}_{i\in[n]}$ and $\{\nu_j\}_{j\in[m]}$ are strictly positive if $\phipi=+\infty$, otherwise, if $\phipi<+\infty$, they can be equal to zero.
The formulation presented in \eqref{eq:matrix-problem}, allows us to express the optimality conditions as follows.

\begin{theorem}
\label{thm:kkt}
A transportation plan $\pi\geq0$ is optimal in \eqref{eq:matrix-problem} if and only if there exist $a\in\R^n$ and $b\in\R^m$ such that
\begin{equation}
 a_i\in\partial \varphi_{\mu_i}((\Pruno \pi)_i),\qquad
 b_j\in\partial \varphi_{\nu_j}((\Prdue\pi)_j)
 \label{eq:potentials}
\end{equation}
and, for every $i\in[n],j\in[m]$, it holds
\begin{equation}
 c_{ij}+a_i+b_j\geq0,\qquad
 \pi_{ij}(c_{ij}+a_i+b_j)=0.
 \label{eq:complementarity}
\end{equation}
\end{theorem}

\begin{proof}
Put $C=\R_+^{n\times m}$ and $J(\pi)=\pair{c}{\pi}+H(L\pi)+\ind_C(\pi)$, where
\begin{equation*}
 L\pi=(\Pruno\pi,\Prdue\pi),\qquad
 H(r,s)=\sum_{i\in[n]} \varphi_{\mu_i}(r_i)+\sum_{j\in[m]} \varphi_{\nu_j}(s_j).
\end{equation*}
For $\pi^0=\one_n\one_m^\top$,
\begin{equation*}
 \pi^0\in\operatorname{int}C,\qquad
 L\pi^0=(m\one_n,n\one_m)\in\operatorname{int}(\dom H).
\end{equation*}
Thus $H$ is continuous at $L\pi^0$ and $H\circ L$ is continuous at
$\pi^0\in C$. The convex chain and sum rules give
\begin{align*}
 \partial J(\pi)&=c+L^\star\partial H(L\pi)+N_C(\pi),\\
 \partial H(r,s)&=
   \prod_i\partial \varphi_{\mu_i}(r_i)\times
   \prod_j\partial \varphi_{\nu_j}(s_j),\\
 (L^\star(a,b))_{ij}&=a_i+b_j,\qquad\qquad
 N_C(\pi)=\{z:z_{ij}\leq0,\ z_{ij}\pi_{ij}=0\text{ for all }i,j\}.
\end{align*}
Consequently,
\begin{align*}
 \pi\in\argmin J
 &\Leftrightarrow 0\in\partial J(\pi)\\
 &\Leftrightarrow
 \left\{\begin{array}{l}
 \text{there exist }a_i\in\partial \varphi_{\mu_i}((\Pruno\pi)_i),\quad
 b_j\in\partial \varphi_{\nu_j}((\Prdue\pi)_j),\\
 -(c_{ij}+a_i+b_j)\leq0,\quad
 \pi_{ij}(c_{ij}+a_i+b_j)=0\quad\text{for all }i,j.
 \end{array}\right.
\end{align*}
In particular, these conditions imply, for every $\gamma\geq0$,
\begin{align*}
 F_c(\gamma)-F_c(\pi)
 &\geq\pair{c}{\gamma-\pi}
      +\pair{a}{\Pruno\gamma-\Pruno\pi}
      +\pair{b}{\Prdue\gamma-\Prdue\pi}\\
 &=\sum_{i,j}(c_{ij}+a_i+b_j)(\gamma_{ij}-\pi_{ij})\\
 &=\sum_{i,j}(c_{ij}+a_i+b_j)\gamma_{ij}\geq0.
\end{align*}
\end{proof}

\subsection{Solution Entry-wise Bounds and Elimination of Expensive Arcs}
\label{subsec:entry-bounds}

In this section, we provide an upper and lower bound on each entry of a solution to the problem in \eqref{eq:matrix-problem}.
For $q\geq0$, let us define the following auxiliary function
\begin{equation}
 \Tij(q):=\max\argmin_{t\geq0}
       \big\{qt+\varphi_{\mu_i}(t)+\varphi_{\nu_j}(t)\big\} \quad\quad (i\in[n], \ j\in[m]).
 \label{eq:scalar-bound}
\end{equation}
Under our assumptions, the scalar objective in \eqref{eq:scalar-bound} is proper, convex, and lower semicontinuous in $t$.
Moreover, it is coercive since it grows at infinity with strictly positive asymptotic slope $q+2\phipi$.
The $\argmin$ set defining $\Tij$ is then non empty and closed, thus $\Tij(q)$ is always finite and well defined. 
Notice that if $\mu_i=\nu_j=0$, we have $\Tij(q)=0$ for every $q\ge 0$.
\begin{theorem}
\label{prop:entry-bound}
    The function $\Tij$ is nonincreasing. 
    Moreover, every optimizer of
    \eqref{eq:matrix-problem} satisfies
    \begin{equation}
     0\leq\pi_{ij}\leq \Tij(c_{ij})\quad\quad (i\in[n],\ j\in[m]).
     \label{eq:entrywise-bound}
    \end{equation}
\end{theorem}

\begin{proof}
    First, we show that
    \begin{equation}
        \Tij(q)= \sup\big(\{t>0:q+d_i^+(t)+e_j^+(t)\leq0\}\cup\{0\}\big),
        \label{eq:def_of_Tij}
    \end{equation}
    where 
    \begin{equation}
    \label{eq:partial_derivatives}
     d_i^+(t)=
     \begin{cases}\phi^\prime_+(t/\mu_i),&\mu_i>0,\\
     \phipi,&\mu_i=0,
     \end{cases}
     \qquad
     e_j^+(t)=
     \begin{cases}\phi^\prime_+(t/\nu_j),&\nu_j>0,\\
     \phipi,&\nu_j=0.
     \end{cases}
    \end{equation}
    Indeed, for the sake of simplicity, let us set $g_q(t)=qt+\varphi_{\mu_i}(t)+\varphi_{\nu_j}(t)$ and 
    \begin{equation}
        \label{eq:alt_T_def}
        T=\Tij(q)=\max\argmin_{t\ge 0} g_q(t).
    \end{equation}
    For $t>0$, we have $g^\prime_{q,+}(t)=q+d_i^+(t)+e^+_j(t)$.
    {\color{black}Convexity gives}
    \begin{align}
     \nonumber 0<t<T&\quad\Longrightarrow\quad
     g^\prime_{q,+}(t)\leq\frac{g_q(T)-g_q(t)}{T-t}\leq0,\\
     \label{eq:lowerboundT} t>T&\quad\Longrightarrow\quad
     0<\frac{g_q(t)-g_q(T)}{t-T}\leq g^\prime_{q,-}(t)\leq g^\prime_{q,+}(t),
    \end{align}
    hence, if $T>0$, we have
    \begin{equation*}
     \left]0,T\right[\subseteq \{t>0:g^\prime_{q,+}(t)\leq0\}\subseteq\left]0,T\right]
    \end{equation*}
    and $\{t>0:g^\prime_{q,+}(t)\leq0\}=\varnothing$ if $T=0$, proving \eqref{eq:def_of_Tij}.
    To conclude that $\Tij$ is monotone, we notice that, owing to the representation of $\Tij$ in \eqref{eq:def_of_Tij}, if $q_1\le q_2$ then
    \begin{equation*}
         g^\prime_{q_1,+}(t)\leq g^\prime_{q_2,+}(t)\quad(t\ge 0)\quad\Longrightarrow\quad
           T_{ij}(q_2)\leq T_{ij}(q_1).
    \end{equation*}
    We now prove the last part of the statement.
    Let $\pi$ be optimal.
    If $\pi_{ij}=0$, then \eqref{eq:entrywise-bound} trivially holds.
    Otherwise, let us set $s=\pi_{ij}>0$ and define $d_i^-,e_j^-$ as in \eqref{eq:partial_derivatives} by replacing the right derivatives with left derivatives for every $t>0$.
    Then, we have $g^\prime_{q,-}(t)=q+d_i^-(t)+e_j^-(t)$.
    Toward a contradiction, assume that $s>\Tij(c_{ij})$.
    By \eqref{eq:lowerboundT}, we have $0<g^\prime_{c_{ij},-}(s)$.
    From Theorem~\ref{thm:kkt}, we find a couple of potentials $a_i$ and $b_j$ for which it holds
    \begin{align*}
     a_i&\geq d_i^-((\Pruno\pi)_i)\geq d_i^-(\pi_{ij}) = d_i^-(s),\\
     b_j&\geq e_j^-((\Prdue\pi)_j)\geq  e_j^-(\pi_{ij}) = e_j^-(s).
    \end{align*}
    Since $\pi_{ij}>0$, we infer that \begin{equation}
        0=c_{ij}+a_i+b_j \geq c_{ij}+d_i^-(s)+e_j^-(s)=g^\prime_{c_{ij},-}(s),
    \end{equation}
    which contradicts $0<g^\prime_{c_{ij},-}(s)$, hence \eqref{eq:entrywise-bound} holds.
\end{proof}

\begin{remark}[Smooth penalties]
If $\mu_i,\nu_j>0$ and $\phi$ is differentiable, we have that $T_{ij}(q)=t>0$ if and only if
\begin{equation*}
 q+\phi^\prime(t/\mu_i)+\phi^\prime(t/\nu_j)=0,
\end{equation*}
otherwise $T_{ij}(q)=0$.
Notice however that no differentiability is needed in Theorem~\ref{prop:entry-bound}.
For reference, if $\phi(t)=t\log t-t+1$, then $\Tij(q)=\sqrt{\mu_i\nu_j}\exp(-q/2)$.
\end{remark}

Notice that if $\Tij(c_{ij})=0$ for some of the costs $c_{ij}$, Theorem \ref{prop:entry-bound} establishes that $\pi_{ij}=0$ for every solution to \eqref{eq:matrix-problem}, allowing us to remove the variable from the problem.
In the next theorem, we give a characterization of this phenomenon.
In what follows, we set
\begin{equation}
 \minr:=\min_{j\in [m]}c_{ij} \qquad \text{and}\quad 
 \minc:=\min_{i\in [n]}c_{ij}.
 \label{eq:min-costs}
\end{equation}

\begin{theorem}[Support bounds]
\label{thm:cutoffs}
Every optimal transportation plan has $\pi_{ij}=0$ if $\mu_i=\nu_j=0$.
Moreover, the following implications hold for every optimizer
\begin{enumerate}[label=(\roman*)]
 \item \label{option:1} If $\phipz > -\infty$, then
 \begin{equation}
  c_{ij}>-2\phipz\quad\Longrightarrow\quad\pi_{ij}=0.
  \label{eq:deletion-cutoff}
 \end{equation}
 \item \label{option:2} If $\phipi<+\infty$, then
 \begin{equation}
  c_{ij}>\minr+\minc+2\phipi
  \quad\Longrightarrow\quad\pi_{ij}=0.
  \label{eq:rerouting-cutoff}
 \end{equation}
\end{enumerate}
\end{theorem}

\begin{proof}
Without loss of generality, let $\pi_{ij}>0$, then $(\Pruno\pi)_i,(\Prdue\pi)_j>0$.
Owing to Theorem \ref{thm:kkt}, we find two potentials, namely $a_i$ and $b_j$, such that
\begin{equation*}
       c_{ij}+a_i+b_j=0 \quad \Rightarrow \qquad c_{ij}=-a_i-b_j.
\end{equation*}
If $\mu_i=\nu_j=0$, the result follows from Theorem \ref{thm:solution_support_characterization}.
Let us now consider \ref{option:1}.
Owing to Theorem \ref{thm:kkt}, we have that $a_i\in\partial\varphi_{\mu_i}((\Pruno\pi)_i)$ and $b_j\in\partial\varphi_{\nu_j}((\Prdue\pi)_j)$.
If $\pi_{ij}>0$, combining \eqref{eq:derivative-order} and
\eqref{eq:marginal-subdifferential} with $c_{ij}+a_i+b_j=0$ gives
\begin{equation*}
 a_i\geq \phipz,\qquad b_j\geq \phipz,
 \qquad c_{ij}=-a_i-b_j\leq-2\phipz
\end{equation*}
which concludes the proof of \ref{option:1}.

Let us now consider \ref{option:2}.
Owing again to Theorem \ref{thm:kkt}, we have $a_l\in\partial\varphi_{\mu_l}((\Pruno\pi)_l)$, thus $a_l\le\phipi$ for every $l\in[n]$ and, likewise, $b_k\leq\phipi$ for every $k\in[m]$.
Choose $j_i\in[m]$ such that $c_{i j_i}=\minr$ and $i_j\in[n]$ such that $c_{i_j j}=\minc$.
Then the optimality conditions yield
\begin{align*}
 0\leq\minr+a_i+b_{j_i}&\quad\Longrightarrow\quad
 a_i\geq-\minr-b_{j_i}\geq-\minr-\phipi,\\
 0\leq\minc+a_{i_j}+b_j&\quad\Longrightarrow\quad
 b_j\geq-\minc-a_{i_j}\geq-\minc-\phipi,
\end{align*}
which, combined together gives $c_{ij}=-a_i-b_j\leq\minr+\minc+2\phipi$,
which concludes the proof.
\end{proof}

Under additional assumptions, the bounds presented in Theorem \ref{thm:cutoffs} simplify as follows.

\begin{corollary}
    \label{thm:cutoffs_simple}
    Under the hypothesis of Theorem \ref{thm:cutoffs}, if $\XC=\YC$, $c(x,x)=0$ for all $x\in\XC$, we have
    \begin{equation}
     (i\in[n],j\in[m])\quad\quad c_{ij}>2\min\{-\phipz,\phipi\}\quad\Longrightarrow\quad\pi_{ij}=0.
     \label{eq:combined-cutoff}
    \end{equation}
\end{corollary}

\begin{proof}
    It follows from Theorem \ref{thm:cutoffs} since $\minr=\minc=0$ for every $i\in[n]$ and every $j\in[m]$.
\end{proof}

It is important to notice that the support bounds presented in Corollary \ref{thm:cutoffs_simple} and in Theorem \ref{thm:cutoffs} do not depend on the input measures $\mu$ and $\nu$, but only on $\phipi$, $\phipz$, and $c$.
The second consequence of the bounds presented in Theorem \ref{prop:entry-bound} is that there exists a solution to \eqref{eq:matrix-problem} whose total mass, intended as the sum of all the values $\pi_{ij}$, is bounded.

\begin{proposition}
\label{prop:bounded-optimizer}
For every optimal transport plan $\pi$ for problem \eqref{eq:matrix-problem}, there is an optimal transport plan, namely $\widehat\pi$, such that $\widehat\pi\leq\pi$ componentwise and, for every couple of indices $(i,j)\in[n]\times[m]$, it holds
\begin{equation}
 \widehat\pi_{ij}>0
 \quad\Longrightarrow\quad
 (\Pruno\widehat\pi)_i\leq\mu_i
 \;\; \text{or}\;\; (\Prdue\widehat\pi)_j\leq\nu_j.
 \label{eq:trimmed-property}
\end{equation}
In particular, we have that $\widehat\pi_{ij}\leq\max\{\mu_i,\nu_j\}$, therefore
\begin{equation}
\|\widehat\pi\|_1\leq B_0:=\sum_{i=1}^n\sum_{j=1}^m\max\{\mu_i,\nu_j\}.
 \label{eq:bounded-optimum}
\end{equation}
\end{proposition}

\begin{proof}
    Let $\pi$ be the solution to \eqref{eq:matrix-problem}.
    If $\pi_{ij}=0$ for every $i\in[n],j\in [m]$, we have nothing to prove.
    Assume then that $i\in[n]$ and $j\in[m]$ are such that $\pi_{ij}>0$ and both the following conditions hold
    \begin{equation}
        (\Pruno\pi)_i>\mu_i \quad\quad\text{and}\quad\quad(\Prdue\pi)_j>\nu_j.
    \end{equation}
    Then, set $\varepsilon=\min\{\pi_{ij},(\Pruno\pi)_i-\mu_i,(\Prdue\pi)_j-\nu_j\}>0$.
    We have
    \begin{enumerate*}[label=(\roman*)]
        \item $\pi_{ij}-\varepsilon\ge 0$, and
        \item $(\Pruno\pi)_i-\varepsilon\ge\mu_i$ and $(\Prdue\pi)_j-\varepsilon\ge\nu_j$.
    \end{enumerate*}
    Then, $\pi^\prime=\pi-\varepsilon\delta_{(x_i,y_j)}$ is feasible for problem \eqref{eq:matrix-problem}.
    Moreover, since $c_{ij}\ge 0$ and owing to the assumptions on $\phi$, we have that $\UTphi(\pi^\prime)\le\UTphi(\pi)$, hence $\pi^\prime$ is also optimal.
    By iterating the process, we build a solution to problem \eqref{eq:matrix-problem}, such that either $(\Pruno\pi)_i\le\mu_i$ or $(\Prdue\pi)_j\le\nu_j$ hold.
    Let us assume that $\pi$ is a solution to \eqref{eq:matrix-problem} such that \eqref{eq:trimmed-property} holds, then
    \begin{equation*}
     \pi_{ij}>0\quad\Longrightarrow\quad \pi_{ij}\leq (\Pruno\pi)_i\leq\mu_i
     \quad\text{or}\quad
     \pi_{ij}\leq (\Prdue\pi)_j\leq\nu_j.
    \end{equation*}
    Thus $\pi_{ij}\leq\max\{\mu_i,\nu_j\}$ and hence \eqref{eq:bounded-optimum}.
\end{proof}

In particular, problem \eqref{eq:matrix-problem} can be cast as a minimization problem over a scaled generalized simplex in $\erre^{n\times m}$, indeed
\begin{equation}
    v_c=\min_{\pi\in\Lambda_{B_0}}F_c(\pi)\quad\quad \text{with}\quad\Lambda_{B_0}:=\{\pi\geq0:\|\pi\|_1\leq B_0\}.
 \label{eq:matrix-simplex}
\end{equation}
Notice that the set $\Lambda_{B_0}$ depends only on $\mu$ and $\nu$ and not on $\Dphi$ and $c$.
To conclude, we infer a lower bound on the marginals of any optimal solution.
First, we define the lower generalized inverse
\begin{equation}
 \Iphi(q):=\inf\{t\geq0:\phi'_+(t)\geq q\},
 \qquad q\in\R,
 \label{eq:lower-inverse}
\end{equation}
with $\inf\varnothing=+\infty$. 
Notice that if $q<0$, the set defining $\Iphi(q)$ is always nonempty as $\phi^\prime_+(1)\ge 0$.

\begin{theorem}
\label{thm:lower-marginals}
Let $U_i,V_j$ be finite, strictly positive numbers such that every $\pi\in\Lambda_{B_0}$ solution to \eqref{eq:matrix-simplex} satisfies $(\Pruno\pi)_i\leq U_i$ and $(\Prdue\pi)_j\leq V_j$.\footnote{Notice that one can always take $U_i=V_j=B_0$.}
Define
\begin{equation}
 A_i:=
 \begin{cases}\phi^\prime_+(U_i/\mu_i),&\mu_i>0,\\
 \phipi,&\mu_i=0,
 \end{cases}
 \qquad
 B_j:=
 \begin{cases}\phi^\prime_+(V_j/\nu_j),&\nu_j>0,\\
 \phipi,&\nu_j=0.
 \end{cases}
 \label{eq:potential-upper-bounds}
\end{equation}
Then every optimizer $\pi$ of \eqref{eq:matrix-problem}, that belongs to $\Lambda_{B_0}$, satisfies
\begin{align}
 (\Pruno\pi)_i&\geq\mu_i\max_{j\in[m]}\Iphi(-c_{ij}-B_j)
 &&\text{if }\mu_i>0,\label{eq:row-lower}\\
 (\Prdue\pi)_j&\geq\nu_j\max_{i\in[n]}\Iphi(-c_{ij}-A_i)
 &&\text{if }\nu_j>0.\label{eq:column-lower}
\end{align}
If $\phipz=-\infty$, these lower bounds are strictly positive.
In particular, there exists a common $\varepsilon>0$ for which if $\mu_i>0$ then $(\Pruno\pi)_i\ge \varepsilon$ and if $\nu_j>0$ then $(\Prdue\pi)_j\ge \varepsilon$ for every optimal $\pi\in\Lambda_{B_0}$.
\end{theorem}

\begin{proof}
Fix an arbitrary optimizer $\pi$ in $\Lambda_{B_0}$. 
Write $r=\Pruno\pi$, $s=\Prdue\pi$, and choose potentials $a,b$ for $\pi$ as in Theorem \ref{thm:kkt}.
Since $r_i\leq U_i$ and $s_j\leq V_j$, we have
\begin{equation*}
    a_i\leq A_i,\qquad \text{and} \qquad b_j\leq B_j.
\end{equation*}
Thus, given $i\in[n]$ such that $\mu_i>0$, for every $j\in[m]$, it holds
\begin{equation*}
 \phi^\prime_+(r_i/\mu_i) \geq a_i\geq-c_{ij}-b_j\geq-c_{ij}-B_j,
\end{equation*}
and hence
\begin{equation*}
 r_i/\mu_i \in\{t\geq0:\phi^\prime_+(t)\geq-c_{ij}-B_j\},\qquad \Rightarrow\qquad r_i\geq\mu_i \Iphi(-c_{ij}-B_j).
\end{equation*}
so that, by taking the maximum with respect to $j\in[m]$ allows us to conclude \eqref{eq:row-lower}.
Likewise, given $j\in[m]$ such that $\nu_j>0$, for every $i\in[n]$, it holds
\begin{equation*}
 \phi^\prime_+(s_j/\nu_j)
 \geq b_j\geq-c_{ij}-a_i\geq-c_{ij}-A_i,
 \qquad
 \Rightarrow\qquad s_j\geq\nu_j \Iphi(-c_{ij}-A_i),
\end{equation*}
thus, by taking the maximum with respect to $i\in[n]$, we infer \eqref{eq:column-lower}.
If $\phipz=-\infty$, then for each finite $q$,
\begin{equation*}
 \lim_{t\downarrow0}\phi^\prime_+(t)=-\infty
 \quad\Longrightarrow\quad
 (\exists\varepsilon_q>0)(\forall t\in\left[0,\varepsilon_q\right[)\quad
 \phi^\prime_+(t)<q
 \quad\Longrightarrow\quad \Iphi(q)\geq\varepsilon_q>0.
\end{equation*}
Since the $A_i,B_j,c_{ij}$ are finite, all the asserted lower bounds are strictly positive in this case.
The last statement follows from setting $\varepsilon=\varepsilon_{\bar q}\min_{i\in I_+,j\in J_+}\{\mu_i,\nu_j\}$ with 
\[
\bar q=\min_{i,j}\{-c_{ij}-A_i,-c_{ij}-B_j\},
\]
and $I_+=\{i\in[n]\ : \mu_i>0\}$ and $J_+=\{j\in[m]\ : \ \nu_j>0\}$.
\end{proof}

\begin{remark}
    Notice that Theorem \ref{thm:lower-marginals} holds for every set of transportation plans $\pi$ such that $\|\pi\|_1\le KB_0$ with $K$ positive constant.
\end{remark}

\section{Cost truncation and a compressed formulation}
\label{sec:truncation}

We now focus on UOT problems induced by truncated costs, i.e., costs that are capped at a given threshold $R>0$.
This class of problems simplifies the structure of the transportation plan, as we do not need to keep track of how the mass is moved between points whose distance is equal to $R$ \cite{Pele2009}.
The problem induced by the truncated cost has a finite-valued minimizer by the same argument used for the original problem as in \eqref{eq:UOTphi}.
Moreover, owing to Theorem \ref{thm:solution_support_characterization}, there exists at least one solution supported over a discrete set.
Given a cost function $c$ and a positive threshold $R>0$, we then set
\begin{equation}
 \ctrun_{ij}:=\min\{c_{ij},R\},\qquad
 \ER:=\{(i,j):c_{ij}<R\},\qquad
 \HR:=([n]\times[m])\setminus \ER,
 \label{eq:capped-cost}
\end{equation}
and define
\begin{equation}
 \hat v_R:=\min_{\pi\geq0}\UTRphint(\pi).
 \label{eq:capped-problem}
\end{equation}
Notice that
\begin{equation}
 R_1\leq R_2\quad\Longrightarrow\quad c^{(R_1)}_{ij}\le c^{(R_2)}_{ij} \quad (i\in[n],\ j\in[m]) \quad\Longrightarrow\quad
 \hat v_{R_1}\leq \hat v_{R_2}\leq v_c.
 \label{eq:value-order}
\end{equation}
The grids selected for $c$ also suffice for its continuous capped cost $\min\{c,R\}$: every selected minimizer of $c(x,\cdot)$ or $c(\cdot,y)$ is still a minimizer after capping.
The problem induced by the truncated cost does not remove the variables in $\HR$, but makes them all have the same cost. 
The following formulation exploits this fact by considering the row and column totals associated with the couples in $\HR$ instead of their individual entries.

\subsection{Compression of the capped links}
\label{subsec:compression}

Let $\eta\in\R_+^{n\times m}$ be supported on $\ER$, that is $\eta_{ij}=0$ for every $(i,j)\in\HR$ and let $u\in\R_+^n$ and $w\in\R_+^m$ satisfy $\|u\|_1=\|w\|_1$.
We define 
\begin{equation}
    \GR(\eta,u,w):=\pair{c}{\eta}+R\|u\|_1
   +\sum_{i\in[n]} \varphi_{\mu_i}((\Pruno\eta)_i+u_i)
   +\sum_{j\in[m]} \varphi_{\nu_j}((\Prdue\eta)_j+w_j).
 \label{eq:compressed-objective}
\end{equation}
and the following feasible set
\begin{equation}
 \FeasR:=\{(\eta,u,w)\in \R_+^{n\times m}\times \R_+^n \times \R_+^m  : \ \eta_{ij}=0\text{ on }\HR, \ \text{and} \ \|u\|_1=\|w\|_1\}.
 \label{eq:compressed-set}
\end{equation}
Notice that only the entries of $\eta$ in $\ER$ need to be stored, so that the formulation uses $|\ER|+n+m$ nonnegative variables and one balance equation.

For any full transportation plan $\pi\geq0$, we define a canonical element $\Csf(\pi)=(\Ceta(\pi),\Cu(\pi),\Cw(\pi))\in\FeasR$ as follows
\begin{equation}
 \Ceta(\pi)=\eta_{ij}=\begin{cases}\pi_{ij},&(i,j)\in \ER,\\0,&(i,j)\in \HR,\end{cases}
 \label{eq:canonical-compression}
\end{equation}
and 
\begin{equation}
    \Cu(\pi)=u_i=\sum_{j:(i,j)\in \HR}\pi_{ij}, \qquad \Cw(\pi)=w_j=\sum_{i:(i,j)\in \HR}\pi_{ij}.
\end{equation}
It is easy to see that $\UTRphint(\pi)=\GR(\Csf(\pi))$.
Conversely, given an element of $(\eta,u,w)\in\FeasR$, we define
\begin{equation}
 \Lift(\eta,u,w) = \eta + \varrho \quad \text{with arbitrary}\quad \varrho\in\Pi(u,w)=\{\gamma\in\R^{n\times m}_+\ : \ \Pruno\gamma=u,\ \Prdue\gamma=w\},
 \label{eq:residual-couplings}
\end{equation}
notice that $\Pi(u,w)\neq\emptyset$, hence $\Lift$ is well-defined.
However, in general, it holds 
\begin{equation}
    \label{eq:not_equal_lift}
    \GR(\eta,u,w)\ge \UTRphint(\Lift(\eta,u,w)).
\end{equation}
In the next lemma, we show that if $(\eta,u,w)$ is optimal, then \eqref{eq:not_equal_lift} holds as equality regardless of how $\varrho$ is selected in \eqref{eq:residual-couplings}.

\begin{lemma}
\label{lem:residual-separation}
Let $(\eta^\star,u^\star,w^\star)\in\FeasR$ be a minimizer of $G_R$.
Then
\begin{equation}
 (i,j)\in \ER\quad\Longrightarrow\quad u^\star_iw^\star_j=0 \qquad \text{and} \qquad \GR(\eta^\star,u^\star,w^\star)=\UTRphint(\Lift(\eta^\star,u^\star,w^\star)).
 \label{eq:residual-separation}
\end{equation}
\end{lemma}

\begin{proof}
Let $(\eta^\star,u^\star,w^\star)\in\FeasR$ be a minimizer of $G_R$.
If $c_{ij}<R$ and $u^\star_iw^\star_j>0$, define
\begin{equation*}
 \delta=\min\{u^\star_i,w^\star_j\}>0,\quad
 \eta^\prime=\eta^\star+\delta e_i e_j^\top,\quad
 u^\prime=u^\star-\delta e_i,\quad w^\prime=w^\star-\delta e_j,
\end{equation*}
where $e_i$ denotes the $i$-th vector of the Euclidean canonical base.
By construction, $\|u^\prime\|_1=\|w^\prime\|_1$, hence $(\eta^\prime,u^\prime,w^\prime)\in\FeasR$ and
\begin{align*}
 \Pruno\eta^\prime+u^\prime &=\Pruno\eta^\star+u^\star,\quad\qquad\quad \Prdue\eta^\prime+w^\prime=\Prdue\eta^\star+w^\star,\\
 \GR(\eta^\prime,&u^\prime,w^\prime)-\GR(\eta^\star,u^\star,w^\star)
 =c_{ij}\delta-R\delta<0,
\end{align*}
contradicting optimality, therefore $u^\star_iw^\star_j=0$ on $\ER$.
Lastly, given $\varrho\in\Pi(u^\star,w^\star)$ we have
\begin{equation*}
 0\leq \varrho_{ij}\leq\min\{u^\star_i,w^\star_j\}=0\qquad((i,j)\in \ER) \quad\Rightarrow\quad \sum_{i,j}\ctrun_{ij}\varrho_{i,j}=R\sum_{i,j}\varrho_{i,j}=R\sum_{i}u^\star_{i},
\end{equation*}
proving the second identity in \eqref{eq:residual-separation}.
\end{proof}

Owing to the latter result, we are able to prove that the UOT problem induced by the truncated cost at threshold $R>0$ is the same as minimizing $\GR$ over the set $\FeasR$.

\begin{theorem}
\label{thm:compression}
Given two discrete measures $\mu$ and $\nu$, we have that
\begin{equation}
 \min_{\pi\in\R^{n\times m}_+}F_{\ctrun}(\pi)=\min_{(\eta,u,w)\in\FeasR}\GR(\eta,u,w).
 \label{eq:compression-equality}
\end{equation}
Moreover, the following two statements hold
\begin{enumerate}
 \item Given $\pi^\star\in\argmin F_{\ctrun}$, then $\Csf(\pi^\star)\in\argmin \GR$ and
 \item if $(\eta^\star,u^\star,w^\star)\in\argmin \GR$, then $\Lift(\eta^\star,u^\star,w^\star)\in \argmin F_{c^{(R)}}$ regardless of how $\varrho\in\Pi(u^\star,w^\star)$, is chosen.
\end{enumerate}
\end{theorem}

\begin{proof}
    Put $\widetilde v_R=\inf_{\FeasR}G_R$.
    Given any $\pi\in\R^{n\times m}_+$, let $(\eta,u,w)=\Csf(\pi)$.
    Then, by definition $\|u\|_1=\|w\|_1$ and
    \begin{equation}
        (\Pruno\eta)_i+u_i = \sum_{j:(i,j)\in \ER}\pi_{ij} + \sum_{j:(i,j)\in \HR}\pi_{ij}=(\Pruno\pi)_i,
    \end{equation}
    likewise, we infer $(\Prdue\eta)_j+w_j=(\Prdue\pi)_j$, thus
    \begin{equation}
        \label{eq:implication}
        \sum_{\ER}c_{ij}\pi_{ij}+\sum_{\HR}R\pi_{ij} =\scalarp{\ctrun,\pi} \qquad \Rightarrow \qquad  G_R(\Csf(\pi))=F_{\ctrun}(\pi),
    \end{equation}
    so that $\widetilde v_R\leq v_R$.
    Let us now take any $(\eta,u,w)\in\FeasR$ and set $\pi=\Lift(\eta,u,w)$.
    Then, by definition of $\Lift$, we have $\Pruno\pi =\Pruno\eta +u$ and $\Prdue\pi=\Prdue\eta +w$.
    Therefore
    \begin{align}
        F_{\ctrun}(\Lift(\eta,u,w))-G_R(\eta,u,w) &= \scalarp{\ctrun,\varrho}-R\tau = \sum_{i,j}(\ctrun_{ij}-R)\varrho_{ij}\leq 0,
     \label{eq:compression-inequality}
    \end{align}
    so that $v_R\leq\widetilde v_R$ hence \eqref{eq:compression-equality}.
    To conclude, we need to prove that $\Csf$ sends a minimizer of $F_{\ctrun}$ into a minimizer of $\GR$ and, vice-versa, $\Lift$ sends a minimizer of $\GR$ in a minimizer of $F_{\ctrun}$ for every possible choice of $\varrho\in\Pi(u,w)$.
    The first implication comes from the identity on the right hand side of \eqref{eq:implication} and from $\tilde v_R=v_R$.
    For the second implication, let us take $(\eta^\star,u^\star,w^\star)$ such that $G_R(\eta^\star,u^\star,w^\star)=\widetilde v_R$.
    Owing to \eqref{eq:compression-inequality}, we have
    \begin{equation*}
     v_R\leq F_{\ctrun}(\Lift(\eta^\star,u^\star,w^\star))\leq G_R(\eta^\star,u^\star,w^\star)=\widetilde v_R=v_R,
    \end{equation*}
    for every $\varrho\in\Pi(u^\star,w^\star)$, hence $\Lift(\eta^\star,u^\star,w^\star)$ is a minimizer for $F_{\ctrun}$.
    \end{proof}

    Notice that Theorem \ref{thm:compression} allows us to identify the set of minimizers of the two problems.
    As a straightforward consequence, we infer the following.

    \begin{corollary}
    \label{cor:bounded-compression}
    Denoted by
    \begin{equation}
     \setKR:=\{(\eta,u,w)\in\FeasR \ : \  \|\eta\|_1\leq B_0,\ \|u\|_1\leq B_0\},
     \label{eq:bounded-compressed-polytope}
    \end{equation}
    we have $v_R=\min_{(\eta,u,w)\in\setKR}G_R(\eta,u,w)$.
    \end{corollary}

    \begin{proof}
    By Proposition~\ref{prop:bounded-optimizer}, there exists a $\pi\in\R^{n\times m}_+$ such that $F_{\ctrun}(\pi)=v_R$ and $\|\pi\|_1\leq B_0$.
    Set $(\eta,u,w)=\Csf(\pi)$.
    Owing to Theorem \ref{thm:compression}, $\Csf(\pi)$ is optimal and $(\eta,u,w)\in\setKR$ since
    \begin{equation*}
        \|\eta\|_1+\tau =\sum_{\ER}\pi_{ij}+\sum_{\HR}\pi_{ij} =\|\pi\|_1\leq B_0,
    \end{equation*}
    which allows us to conclude the proof.
    \end{proof}

\subsection{Error Bounds}
\label{subsec:recovery}
In this section, we study the UOT problem induced by a truncated cost function as an approximation term of the untruncated UOT problem.
In particular, the final goal of the section is to characterize the error induced by approximating $v_c$ with $v_R$.

\subsubsection{Recovering plans for the original cost}

As we saw in Section \ref{subsec:compression}, given an element $(\eta,u,w)\in\FeasR$.
The function $\Lift$ assigns every element of $\FeasR$ to a feasible plan $\pi$ for the original problem and maps a minimizer of $\GR$ to a minimizer of $F_{\ctrun}$.
We describe two additional recovery procedures that define a feasible point for the original UOT problem as in \eqref{eq:matrix-problem} starting from an element of $\FeasR$.
This will allow us to estimate the gap between $v_c$ and $v_R$ via three different error bounds.

\paragraph{Deleting the residuals.}
Given $(\eta,u,w)\in\FeasR$, define
\begin{equation}
 \Dsf(\eta,u,w):=\eta.
 \label{eq:delete-recovery}
\end{equation}
So that $\Dsf$ builds a feasible plan for the original transportation plan by discarding the residual variables $(u,w)$ and keeping only the mass transported at distance smaller than $R$.
Notice that if $\phipz>-\infty$, we must have $\phi(0)<+\infty$, $\Dsf(\eta,u,w)$ has finite marginal penalties, thus $\UTphi(\Dsf(\eta,u,w))<+\infty$. 
Otherwise, $\UTphi(\Dsf(\eta,u,w))$ could attain unbounded value.

\paragraph{Redirecting the residuals.}
Given $(\eta,u,w)\in\FeasR$, define
\begin{equation}
    \Bsf(\eta,u,w)=\eta +Q^u+Q^w, \qquad Q^u := \sum_i u_i\, e_i e_{j_i}^{\top},\qquad \text{and} \quad Q^w := \sum_j w_j\, e_{i_j} e_j^{\top}
    \label{eq:redirect}
\end{equation}
where $j_i \in \operatorname*{argmin}_{j\in[m]} c_{ij}$ and $i_j \in \operatorname*{argmin}_{i\in[n]} c_{ij}$.
Since both $Q^u$ and $Q^w$ are positive matrices, $\Bsf(\eta,u,w)$ is always feasible for \eqref{eq:matrix-problem}.
Unlike $\Dsf$, this construction retains all the residuals and maps it to the closest point available to the source. 
Notice, however, that the marginals of $\Bsf$ are different from the ones of $\Lift$, indeed
\begin{equation}
 \Pruno\Bsf(\eta,u,w)=\Pruno\eta+u+\Pruno Q^w,\qquad
 \Prdue\Bsf(\eta,u,w)=\Prdue\eta+w+\Prdue Q^u.
 \label{eq:redirect-marginals}
\end{equation}

\subsubsection{Constructive Bounds}

Having introduced our three reconstruction routines $\Lift,\Bsf$, and $\Dsf$, we propose our estimation between $v_c$ and $v_R$.
First, we propose an a posteriori bound, which allow to quantify the gap between $v_R$ and $v_{c,R}$ through the solution to the problem on the right hand side of \eqref{eq:compression-equality}.

\begin{theorem}
\label{thm:error-certificates}
Let $(\eta,u,w)\in\FeasR$ be such that $v_R=\GR(\eta,u,w)$. 
Put $c_{\max}=\max_{i,j}c_{ij}$ and set $\tau=\|u\|_1$.
Let us now consider
\begin{align}
 \mathcal{R}_{\Lift}(\eta,u,w)&:=(c_{\max}-R)_+\tau,
 \label{eq:couple-certificate}\\
 \mathcal{R}_{\Dsf}(\eta,u,w)&:=(-2\phipz-R)_+\tau
 &&\text{if }\phipz> -\infty,\label{eq:delete-certificate}\\
  \mathcal{R}_{\Bsf}(\eta,u,w)&:=
 \left[\sum_i(\minr+\phipi)u_i
       +\sum_j(\minc+\phipi)w_j-R\tau\right]_+
 &&\text{if }\phipi<+\infty\label{eq:redirect-certificate}
\end{align}
where $\minc$ is the vector containing the minimum entry per each column and $\minr$ is the vector containing the minimum element for each row.
Then, for each function $\Xsf\in\{\Lift,\Bsf,\Dsf\}$, we have
\begin{equation}
 0\leq F_c(\Xsf(\eta,u,w))-v_c \leq F_c(\Xsf(\eta,u,w))-v_R\leq \mathcal{R}_{\Xsf}(\eta,u,w).
 \label{eq:primal-certificate}
\end{equation}
In particular, for any optimal $(\eta,u,w)$, we have
\begin{equation}
 0\leq v_c-v_R\leq \min\{\mathcal{R}_\Xsf(\eta,u,w)\ : \ \Xsf\in\{\Lift,\Bsf,\Dsf\}\}.
 \label{eq:value-certificate}
\end{equation}
\end{theorem}

\begin{proof}
Let $(\eta,u,w)$ be such that $\GR(\eta,u,w)=v_R$, we then set $r=\Pruno\eta$ and $s=\Prdue\eta$. 
For every finite-valued feasible $\gamma$ for problem \eqref{eq:matrix-problem}, we have
\begin{equation}
 0\leq F_c(\gamma)-v_c\leq F_c(\gamma)-v_R \qquad \text{and} \qquad 0\leq v_c-v_R\leq F_c(\gamma)-v_R.
 \label{eq:recovery-value-chain}
\end{equation}
We now prove the bound in \eqref{eq:primal-certificate} one function $\Xsf\in\{\Lift,\Bsf,\Dsf\}$ at the time.
The statement will then follow by taking the minimum over the three valid bounds.

\emph{Case one: $\Xsf=\Lift$.} 
Let $\varrho\in\Pi(u,w)$, such that $\Lift(\eta,u,w)=\eta+\varrho$.
Then, Lemma \ref{lem:residual-separation} gives $\varrho|_{\ER}=0$ and $F_{\ctrun}(\eta+\varrho)=v_R$.
Hence
\begin{align}
 F_c(\eta+\varrho)-v_R
 &= F_c(\eta+\varrho)-F_{\ctrun}(\eta+\varrho)\notag\\
 &=\sum_{\HR}(c_{ij}-R)\varrho_{ij}
 \leq(c_{\max}-R)_+\sum_{\HR}\varrho_{ij} =(c_{\max}-R)_+\tau.
 \label{eq:exact-coupling-gap}
\end{align}
Notice that \eqref{eq:exact-coupling-gap} implies that if $\tau=0$, then $v_c=v_R$ and that 
\begin{equation*}
    \min\{\mathcal{R}_\Xsf(\eta,u,w)\ : \ \Xsf\in\{\Lift,\Bsf,\Dsf\}\}=0,
\end{equation*}
so that in what follows, we will assume $\tau>0$.

\emph{Case two: $\Xsf=\Bsf$.} 
Assume that $\phipi<+\infty$ as if $\phipi=+\infty$, we have $\mathcal{R}_{\Bsf}<+\infty$ only if $\tau=0$, in which case the bound has been proven already.
Let us set $\rho=\Pruno Q^w$ and $\sigma=\Prdue Q^u$.
Then, owing to \eqref{eq:redirect-marginals}
\begin{align*}
 \varphi_{\mu_i}(r_i+u_i+\rho_i)-\varphi_{\mu_i}(r_i+u_i)&\leq\phipi\rho_i,\quad\text{and}\quad \varphi_{\nu_j}(s_j+w_j+\sigma_j)-\varphi_{\nu_j}(s_j+w_j)&\leq\phipi\sigma_j.
\end{align*}
The latter inequalities imply  
\begin{align*}
 F_c(\Bsf(\eta,u,w))-v_R & =\sum_i\minr u_i+\sum_j\minc w_j-R\tau\\
 &\quad+\sum_i\big[\varphi_{\mu_i}(r_i+u_i+\rho_i)-\varphi_{\mu_i}(r_i+u_i)\big]\\
 &\quad+\sum_j\big[\varphi_{\nu_j}(s_j+w_j+\sigma_j)-\varphi_{\nu_j}(s_j+w_j)\big]\\
 &\leq\sum_i\minr u_i+\sum_j\minc w_j-R\tau + \phipi\sum_i\rho_i+\phipi\sum_j\sigma_j\\
 &=\sum_i(\minr+\phipi)u_i+\sum_j(\minc+\phipi)w_j-R\tau
 \leq \mathcal{R}_{\Bsf}(\eta,u,w).
\end{align*}

\emph{Case three: $\Xsf=\Dsf$.} 
Suppose $\phipz> -\infty$ as otherwise the bound trivially holds or $\tau=0$.
We recall that for any $m>0$, $x\geq0$, and $t>0$, the convexity of $\phi$ gives
%
\begin{equation*}
 \frac{\varphi_m(x+t)-\varphi_m(x)}{t} = \frac{\phi((x+t)/m)-\phi(x/m)}{t/m} \geq \phipz.
\end{equation*}
Since $\varphi_0$ is linear when $\phipi<+\infty$ the quotient is $\phipi$.
Thus, since $\phipz\le\phipi$ for every $m\ge 0$ we have
\begin{equation}
 \varphi_m(x)-\varphi_m(x+t)\leq-\phipz t\qquad(x,t\geq0).
 \label{eq:deletion-increment}
\end{equation}
It follows that
\begin{align*}
 F_c(\eta)-v_R &=-R\tau+\sum_i\big[\varphi_{\mu_i}(r_i)-\varphi_{\mu_i}(r_i+u_i)\big]
             +\sum_j\big[\varphi_{\nu_j}(s_j)-\varphi_{\nu_j}(s_j+w_j)\big]\\
        &\leq-R\tau-\phipz\sum_i u_i-\phipz\sum_j w_j =(-2\phipz-R)\tau\leq(-2\phipz-R)_+\tau,
\end{align*}
which concludes the proof.
\end{proof}

Notice that the bounds presented in Theorem \ref{thm:error-certificates} are a posteriori, meaning that the guarantees is quantifiable only after we have access to the optimal solution $(\eta,u,w)$.
In the remainder of the section, we exploit the properties of the problem to make the bounds a priori, i.e., computable without having access to an optimal solution to the truncated UOT problem.
As a direct consequence of the latter Theorem, we infer the following conditions under which the optimal value of the two problems coincide $v_c=v_R$ and the minima of the truncated problem induce a minimizer for \eqref{eq:matrix-problem}.

\begin{corollary}[Exact recovery]
\label{cor:exact-recovery}
If $(\eta,u,w)$ is optimal for the truncated problem and $\|u\|_1=0$, then $\pi=\eta$ is a solution to \eqref{eq:matrix-problem} and $v_R=v_c$.
Moreover, we have
\begin{enumerate}[label=(\roman*)]
 \item if $R> c_{\max}$, then $\pi=\eta$ is a solution to \eqref{eq:matrix-problem} and $v_R=v_c$;
 \item if $\phipz> -\infty$ and $R\geq-2\phipz$, then $\pi=\eta$ is a solution to \eqref{eq:matrix-problem} and $v_R=v_c$; and
 \item if $\phipi<+\infty$ and $R\geq\max_i\minr+\max_j\minc+2\phipi$, then $\pi=\eta +Q^u+Q^w$ is a solution to \eqref{eq:matrix-problem} and $v_R=v_c$. If $\minc=\minr=0$ for every $i\in[n]$ and every $j\in[m]$ we have $R\geq2\phipi$.
\end{enumerate}
\end{corollary}

\begin{proof}
The first part of the statement trivially follows from the fact that $\mathcal{R}_{\Lift}(\eta,u,w)=0$ whenever $u$ and $w$ are both null vectors.
For the other assertions, we have that for every $(\eta,u,w)$ it holds
\begin{align*}
 R\geq c_{\max}&\quad\Longrightarrow\quad
 \mathcal{R}_{\Lift}(\eta,u,w)=0, \qquad \qquad \qquad R\geq-2\phipz\quad\Longrightarrow\quad
 \mathcal{R}_{\Dsf}(\eta,u,w)=0,\\
 & R\ge \max_i\minr+\max_j\minc+2\phipi \quad\Longrightarrow\quad
 \mathcal{R}_{\Bsf}(\eta,u,w)=0
\end{align*}
since
\begin{equation*}
    \sum_i(\minr+\phipi)u_i+\sum_j(\minc+\phipi)w_j-R\tau \leq(\max_i\minr+\max_j\minc+2\phipi-R)_+\tau.
\end{equation*}
\end{proof}

We now show that even when the hypothesis of Corollary \ref{cor:exact-recovery} do not hold, we are able to give an upper bound on the error between the truncated and untruncated problems.
First, we introduce a Lemma that bounds the total mass of the minimizer for the truncated problem.

\begin{lemma}
\label{lem:capped-mass}
Let $\pi$ be a transportation plan that minimizes $F_{\ctrun}$, then $\pi_{ij}\leq \Tij(R)$ for every $(i,j)\in \HR$, where $\Tij$ is defined as in \eqref{eq:scalar-bound}.
In particular, for every $R>0$, let $(\eta,u,w)\in\argmin\GR$, then we have
\begin{equation}
 u_i\leq\sum_{j:(i,j)\in \HR}T_{ij}(R),\qquad
 w_j\leq\sum_{i:(i,j)\in \HR}T_{ij}(R),\qquad
 \tau\leq A_R:=\sum_{(i,j) \in \HR}\Tij(R).
 \label{eq:residual-mass-bound}
\end{equation}
\end{lemma}

\begin{proof}
Theorem \ref{prop:entry-bound} gives $\pi_{ij}\leq \Tij(\ctrun_{ij})=\Tij(R)$ for every $(i,j)\in\HR$.
Let $(\eta,u,w)\in\FeasR$ be a minimizer of $\GR$ and let $\pi=\Lift(\eta,u,w)=\eta+\varrho$. 
Then we have $\eta|_{\HR}=0$, hence
\begin{align*}
 u_i =\sum_{j:(i,j)\in \HR}\varrho_{ij} = \sum_{j:(i,j)\in \HR}\pi_{ij} \le \sum_{j:(i,j)\in \HR}T_{ij}(R),  \qquad \tau&=\sum_{(i,j)\in \HR} \pi_{ij}  \leq\sum_{(i,j)\in \HR}T_{ij}(R)=A_R
\end{align*}
a similar argument yields the bound on $w_j$.
\end{proof}

Finally, we present our a priori bounds on the error on the gap between $v_c$ and $v_R$.

\begin{theorem}
\label{cor:apriori}
Define $A_R:=\sum_{(i,j)\in\HR}T_{ij}(R)$, then
\begin{align}
 \mathcal R_{\mathrm{\Lift}}(R)
  &:=\sum_{(i,j)\in \HR}(c_{ij}-R)T_{ij}(R),
  \label{eq:apriori-couple}\\
 \mathcal R_{\Dsf}(R) &:=(-2\phipz-R)_+ A_R&&\text{if }\phipz> -\infty,\label{eq:apriori-delete}\\
 \mathcal R_{\Bsf}(R) &:=\sum_{(i,j)\in \HR}(\minr+\minc+2\phipi-R)_+T_{ij}(R)
  &&\text{if }\phipi<+\infty.\label{eq:apriori-redirect}
\end{align}
In particular, $0\leq v_c-v_R\leq\min\{\mathcal R_\Xsf(R)\ : \ \Xsf\in\{\Lift,\Bsf,\Dsf\}\}$.
\end{theorem}

\begin{proof}
We prove the estimates in \eqref{eq:apriori-couple}, \eqref{eq:apriori-delete}, and \eqref{eq:apriori-redirect} one by one, the last statement of the theorem follows by taking the minimum over the three alternatives.
The first two estimates follow by combining the estimates in \eqref{eq:couple-certificate} and \eqref{eq:delete-certificate} with \eqref{eq:residual-mass-bound}.
We then need to prove only \eqref{eq:apriori-redirect}.
Notice that given an optimal $(\eta,u,w)$ denoted by $\Lift(\eta,u,w)=\eta+\varrho$, we have
\begin{align*}
    \Big[\sum_i(\minr+\phipi)u_i
       +\sum_j(\minc+\phipi)w_j-R\tau\Big]_+&=\Big[\sum_{\HR}(\minr+\minc+2\phipi-R)\varrho_{ij}\Big]_+\\
       &\le \sum_{\HR}(\minr+\minc+2\phipi-R)_+\pi_{ij},
\end{align*}
which, combined with $\pi_{ij}\le \Tij(R)$ for $(i,j)\in\HR$ concludes the proof.
\end{proof}


\section{UOT Induced by Entropy-like Penalties}

In this section, we specialize our discussion to the case in which $\phi$ is one of the powerlike entropies considered in \cite{liero2018optimal}.
Moreover, we assume that $\XC=\YC\subset\R^n$ are discrete and finite and $c(x,y)=h(\|x-y\|)$ with $h:\R_+\to\R_+$ a convex, increasing function such that $h(t)=0$ if and only if $t=0$.
Notice that, owing to Remark \ref{rmk:supp_prod}, the support of the optimal transportation plan is always contained in $(\AC\cup\BC)\times(\AC\cup\BC)$ since $\minr=\minc=0$ and $h(t)>0$ if $t>0$.

\subsection{Powerlike Entropies and Their Properties}

Given $p\in\R$, the powerlike entropy $\Up$ associated with $p$ is defined as
\begin{equation*}
 \Up(t) = 
    \begin{cases}
    t - 1 - \log t &\text{if}\ p=0\\
    t \log t - t + 1 &\text{if}\ p=1\ \\[1ex]
    \dfrac{1}{p(p-1)} (t^p - p(t-1)-1) &\text{otherwise}\ 
    \end{cases}
    \qquad \dom(\Up) = 
    \begin{cases}
    \left[0,+\infty\right[ & \text{if}\ p> 0\\
    \left]0,+\infty\right[ & \text{if}\ p\leq 0
    \end{cases}
\end{equation*}
Moreover, since every $\Up$ is differentiable on $\R_{++}$, we have
\begin{equation*}
     \partial\Up(t) = 
        \begin{cases}
            \{\Up^\prime(t)\} &\text{if}\ t>0\\
        \left]-\infty, \Upo\right] &\text{if}\ t=0\ \text{and}\ p> 1,
        \end{cases}
        \qquad 
        \dom (\partial \Up) = 
        \begin{cases}
            \left[0,+\infty\right[ &\text{if}\ p>1\\
            \left]0,+\infty\right[ &\text{if}\ p\leq 1
        \end{cases}
\end{equation*}
where $\Upo=(\Up^\prime)_+(0)$ is the right derivative of $\Up$ in $0$ and
\begin{equation*}
    (\forall\,t>0)\quad \Up^\prime(t) = 
    \begin{cases}
        \dfrac{1}{p-1}(t^{p-1}-1) & \text{if}\ p\in \R\setminus\{0,1\}\\[1ex]
        \log t & \text{if}\ p=1\\[1ex]
        1- t^{-1} &\text{if}\ p=0
    \end{cases}
\end{equation*}
Note, in addition that
\begin{equation*}
    \Upo=\begin{cases}
        -\infty & \text{if}\ p\le 1\\[1ex]
        \dfrac{1}{1-p} &\text{if}\ p > 1
    \end{cases}
    \quad\quad\quad\text{and}\quad\quad\quad
    \Upi =
    \begin{cases}
        \dfrac{1}{1-p} & \text{if}\ p<1\\[1ex]
        +\infty &\text{if}\ p\geq 1.
    \end{cases}
\end{equation*}

\begin{proposition}\label{prop:15032026}
The powerlike entropy functions $\Up$ satisfy the following properties:
\begin{enumerate}[label=(\roman*)]
    \item For every $0<L<M$, $\Up^\prime$ is Lipschitz continuous on $[L,M]$. Moreover, $(\Up^\prime)_+$ is Lipschitz continuous on
    $[0,M]$ with constant $M^{p-2}$ if $p\ge2$, and $(p-1)$-Hölder continuous with constant $(p-1)^{-1}$ if $1<p<2$.

    \item For $p<2$, $\Up$ is $2$-uniformly convex on $[L,M]$ whenever $0<L<M$. If $p>2$, $\Up$ is $p$-uniformly convex on $[L,\infty)$ for every $L>0$.

    \item $\Udue$ is quadratic, with Lipschitz derivative, and is $2$-uniformly convex on $[0,\infty)$.
\end{enumerate}
In particular, every $\Up$ has Lipschitz derivative and is uniformly convex on each compact interval $[L,M]\subset(0,\infty)$.
\end{proposition}

\begin{proof}
The statements follow directly from the formulas for $\Up^\prime$ and $\Up^{\prime\prime}$.
\end{proof}

\subsection{Truncated UOT with Powerlike Entropies}

Denoting by $\Dp$ the Csiszár divergence induced by $\Up$ and by $\FeasRp$ the set of feasible triplets $(\eta,u,w)$ as in \eqref{eq:compressed-set}, we set
\begin{equation}
\label{eq:uotp}
    v_p:=\min_{\pi_{ij}\ge 0}\UTp(\pi) \quad\quad \text{and} \quad\quad v_{p,R}:=\min_{\FeasRp}\GR(\eta,u,w)
\end{equation}
where 
\begin{align*}
    \UTp(\pi)&=\scalarp{c,\pi}+\Dp(\pi^{(1)}|\mu)+\Dp(\pi^{(2)}|\nu)\\
    \GR(\eta,u,w)&=\scalarp{c,\eta}+R\|u\|_1+\Dp(\Pruno\eta+u|\mu)+\Dp(\Prdue\eta+w|\nu).
\end{align*}
Following the techniques proposed in \cite{beck2017linearly}, it is possible to show that the properties of $\Up$ readily translate to properties over the objective function $\GR$.
Denoted by $\FeasRp^\star$ the set of minimizers of the truncated problem on the right hand side of \eqref{eq:uotp}, we have the following.

\begin{theorem}
\label{lemma:29102025}
    Let $\pi$ be a solution to the UOT problem induced by $\Dp$.
    Then,
    \begin{equation}
        \pi_{ij}\le \Tijp(c_{ij}),\quad\quad\text{where}\quad \Tijp(t):=\begin{cases}
            \frac{\mu_i\nu_j}{(\mu_i^{p-1}+\nu_j^{p-1})^{\frac{1}{p-1}}}\big(2-(p-1)t\big)^{\frac{1}{p-1}}_+ \qquad &\text{if}\quad p > 1,\\[0.3cm]
            \sqrt{\mu_i\nu_j}e^{-\frac{t}{2}}\qquad &\text{if}\quad p= 1,\\[0.3cm]
            \Big(\frac{\mu_i^{1-p}+\nu_j^{1-p}}{2+(1-p)t}\Big)^{\frac{1}{1-p}}\qquad &\text{otherwise}.\\[0.3cm]
        \end{cases}
    \end{equation}
    Moreover, if $p \neq 1$, then $\pi_{ij}=0$ if $c_{ij}>\frac{2}{|1-p|}$.
\end{theorem}

\begin{proof}
    First, we observe that if $\pi_{ij}=0$ all the estimates are trivial, thus we will tacitly assume that $\pi_{ij}>0$.
    Second, we notice that $\Upo=\frac{1}{1-p}$ if $p>1$ and $\Upo=-\infty$ otherwise.
    Likewise $\Upi=\frac{1}{1-p}$ if $p<1$ and $\Upi=\infty$ otherwise.
    Thus, owing to Theorem \ref{thm:cutoffs}, gives
    \begin{enumerate*}[label=(\roman*)]
        \item $\pi_{ij}=0$ whenever $c_{ij}\ge\frac{2}{p-1}$ if $p>1$; and
        \item $\pi_{ij}=0$ if $c_{ij}>\frac{2}{1-p}$ if $p<1$.
    \end{enumerate*} 
    Notice that $\Up$ is differentiable on $\R_{++}$ and $(\Up^\prime)_+(t)=(\Up^\prime)_-(t)=\Up^\prime(t)$ with
    \begin{equation}
        (\forall t\in\R_{++}) \qquad\qquad \Up^\prime(t)=\begin{cases}
            \frac{1}{p-1}\big(t^{p-1}-1\big)\quad\quad &\text{if}\quad p\neq 1 \\
            \log(t)\quad\quad &\text{otherwise.} 
            \end{cases}
    \end{equation}
    Following Theorem \ref{prop:entry-bound}, we compute the function $\Tij$ associated with all the $\Up$.
    First, we consider the case $p>1$.
    Since $(\Up^\prime)_\infty=\infty$ hence, we assume $\mu_i,\nu_j>0$.
    A direct computation shows that $\Tij(t)=s$ if and only if
    \begin{equation}
        -t=\frac{1}{p-1}\bigg(\frac{s^{p-1}}{\mu_i^{p-1}}+\frac{s^{p-1}}{\nu_j^{p-1}}\bigg)-\frac{2}{p-1}
        \quad\iff\quad s^{p-1}=\frac{(\mu_i\nu_j)^{p-1}}{(\mu_i^{p-1}+\nu_j^{p-1})}\big(2-(p-1)t\big),
    \end{equation}
    which allows us to conclude the first part of the proof.
    Let us now consider the case $p<1$.
    Owing to Corollary \ref{thm:cutoffs_simple}, we have that $\pi_{ij}=0$ if $c_{ij} > \frac{2}{1-p}$.
    Let us then consider the case in which $c_{ij}\le \frac{2}{1-p}$.
    In this case, if $\mu_i,\nu_j>0$, we have $\Tij(t)=s$ if and only if
    \begin{equation}
        -t=\frac{2}{1-p}-\frac{1}{1-p}\bigg(\frac{\mu_i^{1-p}+\nu_j^{1-p}}{s^{1-p}}\bigg) \quad\iff\quad 
        s^{1-p}=\frac{\mu_i^{1-p}+\nu_j^{1-p}}{(2-(p-1)t)}.
    \end{equation}
    Via a similar computation we handle the case in which either $\mu_i=0$ or $\nu_j=0$.
    Finally, if $p=1$, we have
    \begin{equation}
        -t=\log(s/\mu_i)+\log(s/\nu_j) \quad\iff\quad s=\sqrt{\mu_i\nu_j}e^{-t/2}.
    \end{equation}
    For all the cases, the proof follows by setting $t=c_{ij}$.
\end{proof}

The bound on the entries of an optimal transportation plan $\pi_{ij}$ can be easily translated in bounds on the error between the objective $v_{p,R}$ and $v_{p}$.

\begin{corollary}
    \label{thm:allcases_p_R_bound}
    Given $p\in\R$, we have that
    \begin{equation}
        v_p-v_{p,R}\le 
        \begin{cases}
            \mathcal{O}\big((\frac{2}{p-1}-R)^{\frac{p}{p-1}}\big)&\text{if}\;p>1\\
            \mathcal{O}\big(e^{-R/2}\big)&\text{if}\;p=1\\
            \mathcal{O}\big(R^{-\frac{1}{1-p}}\big)&\text{otherwise}.
        \end{cases}
    \end{equation}
\end{corollary}

\begin{proof}
    It follows by summing over $\HR$ the bounds presented in Theorem \ref{cor:apriori} combined with the upper bounds on $\Tijp$ presented in Theorem \ref{lemma:29102025}.
\end{proof}

\section{Pairwise Frank-Wolfe for UOT}
\label{sec:PFW}

In this section, we study an efficient implementation of the Pairwise Frank-Wolfe (PFW) algorithm for the UOT problem.
In what follows, we let $P=\operatorname{conv}(\mathcal V)\subset\R^d$ be a non-empty polytope with finite vertex set $\mathcal V$, with cardinality $N$ and denote by $D$ the diameter of $P$. 
Let $F$ be a proper lower semi-continuous convex objective with $\XC^\star:=\operatorname*{argmin}_{x\in P}F(x)\ne\varnothing$ and finite minimum $F^\star:=\min_{x\in P}F(x)$.
Given $x_0\in P$, assume that
\begin{equation}
 \Isubl=\{x\in P:F(x)\le F(x_0)\}
 \label{s6:levelset}
\end{equation}
is contained in an open convex set on which $F$ is continuously differentiable. 
Given $x^{(k)}\in P$, we say that $\xi_{k,v}\ge0$ are the coefficients associated with an explicit representation of $x^{(k)}$ if
\begin{equation}
 x^{(k)}=\sum_{v\in\mathcal V}\xi_{k,v}v,\qquad\sum_v\xi_{k,v}=1,
 \label{s6:representation}
\end{equation}
and denote by $S_k=\{v:\xi_{k,v}>0\}$ the active set of the representation.
Given a point $x^{(k)}$ with explicit representation $\xi_{k,v}$, a PFW iteration is defined by the following two linear minimization oracles (LMOs)
\begin{align}
 s_k&\in\operatorname*{argmin}_{v\in\mathcal V}
                 \langle\nabla F(x^{(k)}),v\rangle,
 &v_k&\in\operatorname*{argmax}_{v\in S_k}
                 \langle\nabla F(x^{(k)}),v\rangle,\nonumber\\
 d_k&=s_k-v_k,
 &g_k&=-\langle\nabla F(x^{(k)}),d_k\rangle,\nonumber\\
 \gamma_k^{\max}&=\xi_{k,v_k},
 &x^{(k+1)}&=x^{(k)}+\gamma_kd_k,
 \label{s6:pfw}
\end{align}
where $0<\gamma_k\le\gamma_k^{\max}$, with $\gamma^{\rm max}_k$ being a feasible maximum stepsize, and $g_k$ is the pairwise Frank-Wolfe gap in $x^{(k)}$.
At each iteration, if $g_k=0$, the algorithm stops, otherwise $s_k\ne v_k$, $x^{(k+1)}\neq x^{(k)}$, and the coefficients of $x^{(k+1)}$ are updated as
\begin{equation}
 \xi_{k+1,s_k}=\xi_{k,s_k}+\gamma_k,\qquad
 \xi_{k+1,v_k}=\xi_{k,v_k}-\gamma_k,
 \label{s6:weights}
\end{equation}
with all other coefficients unchanged and the active set is updated accordingly.

\subsection{The bounded compressed problem and its vertices}

Given $R>0$, let $\setKR$ be the polytope defined in \eqref{eq:bounded-compressed-polytope}, where $B_0$ is the mass bound in \eqref{eq:bounded-optimum}.
We now show that we can factor $\setKR$ as $\setKR=P_\eta\times P_S$ where
\begin{equation}
 \begin{aligned}
 P_\eta&=\{\eta\in\R_+^{n\times m}:\eta_{ij}=0\text{ on }\HR, \ \|\eta\|_1\le B_0\},\\
 P_S&=\{(u,w)\in\R_+^n\times\R_+^m: \ \|u\|_1=\|w\|_1\le B_0\}
 \end{aligned}
 \label{s6:polytope}
\end{equation}
and denote by $z=(\eta,u,w)$ a generic element of $\setKR$.
Notice that only the coordinates of $\eta$ indexed by $\ER$ need to be stored.
Given $z=(\eta,u,w)$, we set
\begin{equation}
 r(z)=\Pruno\eta+u,\qquad s(z)=\Prdue\eta+w,
 \label{s6:marginalmap}
\end{equation}
using the marginal operators from \eqref{eq:matrix-problem}.
Notice that $r(z),s(z)\le 2 B_0$ coordinate-wise for every $z\in\setKR$.
First, we characterize the set of vertices of $\setKR$ and prove that the problems defining the LMOs in \eqref{s6:pfw} can be solved through a simple search when $P=\setKR$.

\begin{proposition}
\label{s6:vertices}
The vertex sets of the two factors in \eqref{s6:polytope} are
\begin{equation*}
 \mathcal V_\eta=\{0\}\cup\{B_0e_{ij}:(i,j)\in \ER\},\qquad
 \mathcal V_S=\{0\}\cup\{B_0(e_i,e_j):i\in[n],\ j\in[m]\}.
\end{equation*}
Consequently $\mathcal V_R=\mathcal V_\eta\times\mathcal V_S$ is the vertex set of $\setKR$. 
In particular, given a linear functional restricted over $\setKR$ with coefficients $(Q,a,b)$, a minimizing vertex for the LMO is obtained as follows:
\begin{enumerate}[label=(\roman*)]
\item Minimize $Q_{ij}$ on $\ER$. Select $B_0e_{ij}$ if its minimum is negative, and select $0$ otherwise. If $\ER=\varnothing$, select $0$.
\item Independently minimize $a_i$ and $b_j$. Select $B_0(e_{\tilde{i}},e_{\tilde{j}})$ if $a_{\tilde i}+b_{\tilde j}=\min_i a_i+\min_j b_j<0$, and select $0$ otherwise.
\end{enumerate}
Ties may be resolved arbitrarily among minimizing vertices.
\end{proposition}

\begin{proof}
The first factor is a scaled simplex with its zero vertex.
For $(u,w)\in P_S$, let $t=\one^\top u=\one^\top w$.
If $t>0$, choose a nonnegative matrix $\rho$ with row sums $u$ and column sums $w$, for example $\rho=uw^\top/t$. 
Then
\begin{equation*}
 (u,w)=\left(1-\frac{t}{B_0}\right)0+
         \sum_{i,j}\frac{\rho_{ij}}{B_0}\,B_0(e_i,e_j).
\end{equation*}
For $t=0$ the point is zero. 
Each displayed nonzero atom is extreme: its zero coordinates force every member of a convex representation to have the same two-coordinate support, and its total mass $B_0$ forces the same mass. 
The zero point is extreme by nonnegativity. 
This proves the second vertex description. 
Vertices of a product are products of vertices.
Linear minimization separates between the factors, and its value at a nonzero residual vertex is $B_0(a_i+b_j)$, proving the oracle formulas.
\end{proof}
Notice that if $P=\setKR$, Proposition~\ref{s6:vertices} allows us to compute $s_k$ without enumerating $\mathcal V_R$.
Moreover, any vertex of $P$ has at most three nonzero coordinates, and at most one new joint atom is added per iteration.

\subsection{Convergence Guarantees for the PFW Algorithm}

In this section, we layout our standing assumptions on the stepsize implemented in the PFW method and on the regularity of the function $F$.
We show that these assumptions ensure the convergence of the PFW and characterize its convergence speed.
First, we state the general set of conditions of $F$ that we will assume throughout the remainder of the section.

\begin{assumption}
\label{assumption_F_P}
    Let $\mathcal P\subset\R^{n}$ be a polytope.
    Given $F\in\Gamma_0(\R^n)$, let $\XC^\star:=\argmin_{z\in \mathcal{P}}F(z)$ be the set of minimizers of $F$ over $\mathcal P$.
    We assume there exists a closed and convex set $\CC$ such that
\begin{enumerate}
    \item $\XC^\star\subset\CC\subset \mathcal P$; 
    \item $F$ is differentiable over an open set containing $\CC$, has an $\alpha$-Hölder continuous gradient on $\CC$ with constant $H$; and
    \item there exists a $\beta\ge 1$, $q\ge2$, and $\kappa_q>0$ for which it holds
    \begin{equation}
    \label{eq:unif_p_cond_assumption}
        \frac{\kappa_q}{q}
        d(z,\XC^\star)^q
        \le
        F(z^\star)-F(z)
        -\beta\scalarp{\nabla F(z), z^\star-z},
    \end{equation}
    where $d(z,\XC^\star)=\min_{z^\star\in\XC^\star}\|z-z^\star\|$ for all $z\in\CC$ and all $z^\star\in\XC^\star$.
\end{enumerate} 
\end{assumption}

\subsubsection{Step sizes and sufficient decrease}

Given $x_0$, let $\Isubl$ be defined as in \eqref{s6:levelset}.
Assume that $\nabla F$ is $\alpha$-Hölder continuous on $\Isubl$, with $0<\alpha\le1$ and constant $H>0$ and set $L=HD^{1+\alpha}$.
In what follows, we consider step sizes satisfying
\begin{equation}
 \gamma_k\ge\bar\gamma_k
 :=\min\{\gamma_k^{\max},\varsigma g_k^{1/\alpha}\},\qquad
 F(x^{(k)})-F(x^{(k+1)})\ge \vartheta\bar\gamma_k g_k.
 \label{s6:stepconditions}
\end{equation}
where $\varsigma,\vartheta\ge 0$ are suitable positive parameters.
We say that a PFW iteration is a good iteration if $\varsigma g_k^{1/\alpha}\le\gamma_k^{\max}$, so that non-good iteration takes the full step $\gamma_k=\gamma_k^{\max}$.
The next proposition ensures that the assumptions on the stepsize are satisfied by three of the most common stepsizes used in iterative routines.
\begin{proposition}
\label{s6:linesearch}
If $\nabla F$ is $\alpha$-Hölder continuous over $\mathcal C$, the following step sizes satisfy \eqref{s6:stepconditions}.
\begin{enumerate}[label=(\roman*)]
\item The model step
\begin{equation*}
 \gamma^{\rm ms}_k=\min\{\gamma_k^{\max},(g_k/L)^{1/\alpha}\},
\end{equation*}
satisfies \eqref{s6:stepconditions} with $\varsigma=L^{-1/\alpha}$ and $\vartheta=\alpha/(1+\alpha)$.
\item The optimal stepsize $\gamma^{\rm opt}_k$ obtained via exact minimization of $F(x^{(k)}+\gamma d_k)$ over $\left[0,\gamma_k^{\max}\right]$ satisfies \eqref{s6:stepconditions} with the same constants as the model stepsize.
\item The Armijo stepsize, denoted by $\gamma_k^{\rm arm}$, is defined by backtracking as follows. Given $0<\rho<1$ and $0<\sigma<1$, start from $\gamma_k^{\rm max}$ and successively multiply $\gamma_k^{\rm max}$ by $\rho$ until $F(x^{(k)}+\gamma d_k)\le F(x^{(k)})-\sigma\gamma g_k$.
Then the Armijo rule satisfies \eqref{s6:stepconditions} with
\begin{equation*}
 \varsigma=\rho\left(\frac{(1+\alpha)(1-\sigma)}L\right)^{1/\alpha},
 \qquad \vartheta=\sigma.
\end{equation*}
\end{enumerate}
\end{proposition}

\begin{proof}
For the sake of simplicity, let us fix $x=x^{(k)}$, $d=d_k$, $g=g_k>0$, and $\gamma^{\max}=\gamma_k^{\max}$.
The descent lemma for functions with $\alpha$-Hölder gradient yields 

\begin{align}
 F(x+td)&\le F(x)-tg+\frac{L}{1+\alpha}t^{1+\alpha},
 \label{s6:descent}
\end{align}
so that for $t_0=\min\{\gamma^{\max},(g/L)^{1/\alpha}\}$ we have
\begin{equation*}
 F(x)-F(x+\gamma^{\rm ms}d) = F(x)-F(x+t_0d)\ge\frac{\alpha}{1+\alpha}t_0g
\end{equation*}
and
\begin{equation}
    \langle\nabla F(x+\gamma d),d\rangle
  \le -g+H\gamma^{\alpha}\|d\|^{1+\alpha}
  \le -g+L\gamma^{\alpha}<0\qquad \text{for} \qquad 0\le\gamma<(g/L)^{1/\alpha},
\end{equation}
which ensures that $h_{k+1}\le h_k$.
Proving the result for $\gamma^{\rm ms}$.
The result on the optimal stepsize follows by the same argument since $F(x)-F(x+\gamma^{\rm opt} d)\ge F(x)-F(x+\gamma^{\rm ms}d)$ and the descent lemma ensures that $\gamma\to F(x+\gamma d)$ is decreasing on $[0,(g/L)^{\frac{1}{\alpha}}]$.
Lastly, we consider the Armijo stepsize, for which the sufficient decrease condition is trivially satisfied by setting $\sigma=\vartheta$.
If the initial trial is accepted, it equals $\gamma^{\max}$, hence both conditions are trivially met.
Owing to the descent lemma, we notice that 
\begin{equation*}
 \gamma_A=\left(\frac{(1+\alpha)(1-\sigma)g}{L}\right)^{1/\alpha} \qquad \Rightarrow \qquad F(x) - F(x+\gamma_A d) \ge\gamma_A\sigma g.
\end{equation*}
In particular, it must be that the Armijo stepsize satisfies $\rho\gamma_A \le \gamma^{\rm arm},$ since if $\gamma^{\rm arm} < \rho \gamma_A$ it means that we would have accepted $\tilde{\gamma}^{\rm arm}=\gamma^{\rm arm}/\rho$ during the backtracking procedure since $\tilde{\gamma}^{\rm arm}<\gamma_A$, contradicting the definition of $\gamma^{\rm arm}$, allowing us to conclude the proof.
\end{proof}

\subsubsection{Convergence Rate Guarantees}
We now characterize the convergence speed of the PFW depending on the parameters $q$ and $\alpha$ defined in Assumption \ref{assumption_F_P}.
For a finite set $V$ and $x\in\operatorname{conv}(V)$, let $\mathcal S_V(x)$ be the sets $S\subset V$ for which $x$ has an explicit representation via the elements in $S$ and the coefficients associated with this representation are all strictly positive coefficients. 
Define
\begin{align*}
 \operatorname{PdirW}(V,r,x)
 &:=\min_{S\in\mathcal S_V(x)}\max_{s\in V,\,v\in S}
                \left\langle\frac r{\|r\|},s-v\right\rangle,\\
 W:=\operatorname{PWidth}(\mathcal V)
 &:=\inf_{\substack{Q\text{ a face of }P,\ x\in Q\\
                    r\in\operatorname{cone}(Q-x)\setminus\{0\}}}
                  \operatorname{PdirW}(\mathcal V\cap Q,r,x).
\end{align*}
This width is strictly positive for a nontrivial polytope, with directions taken in its affine hull. 
{\color{black}We use the following geometric inequality from \cite[Theorem~3]{lacoste2015global}: for every nonoptimal $x\in P$, every active representation of $x$, and every $x^\star\in \XC^\star$, we have
\begin{equation}
 \langle-\nabla F(x),x^\star-x\rangle
 \le \frac{\langle-\nabla F(x),d\rangle}W\|x^\star-x\|,
 \label{s6:width}
\end{equation}
where $d$ is the pairwise Frank-Wolfe direction computed at point $x$.}
We are now ready to state the general convergence result that we will use for the UOT problem.
This extends the results for $2$-conditioned Lipschitz smooth objective functions, given in \cite{lacoste2015global}, to the case in which the objective is $q$-conditioned and Hölder smooth.

\begin{theorem}
\label{thm:general_convergence_PFW}
Let $F$ be such that its gradient is Hölder continuous on $\Isubl$ with exponent $\alpha$ and constant $H$.
Assume that there exist $q\ge2$, $\beta\ge1$, and $\kappa_q>0$ for which \eqref{eq:unif_p_cond_assumption} holds with $\CC=\Isubl$.
Moreover, let $h_k=F(x^{(k)})-F^\star$ be the primal gap, $t_k$ be the number of good iterations among the first $k-1$, and set
\begin{equation}
 A_*:=\frac{W\kappa_q^{1/q}}\beta
          \left(\frac q{q-1}\right)^{(q-1)/q},\qquad
 c_*:=\varsigma \vartheta A_*^{(1+\alpha)/\alpha},\qquad
 s:=\frac{q-1-\alpha}{q\alpha}\ge0.
 \label{s6:rateconstants}
\end{equation}
Then, denoted by $x^{(k)}$ the $k$-th iterate of the PFW algorithm, if $h_0>0$ and the stepsize satisfies \eqref{s6:stepconditions}, we have that
\begin{enumerate}[label=(\roman*)]
\item if $s>0$, then $h_k\le\left(h_0^{-s}+sc_*t_k\right)^{-1/s}$;
\item if $s=0$, equivalently $q=2$ and $\alpha=1$, then $h_k\le h_0(1-\theta)^{t_k}$ with $\theta:=\min\{c_*,1/2\}$.
\end{enumerate}
\end{theorem}

\begin{proof}
Fix a nonoptimal iterate $x$ and choose a nearest optimizer to $x$, namely $x^\star$.
Denoted with $h=F(x)-F^\star$ and $\Delta=\|x-x^\star\|$, conditions
\eqref{eq:unif_p_cond_assumption} and \eqref{s6:width} give
\begin{equation}
 h \le \frac{\beta g}W \Delta-\frac{\kappa_q}{q}\Delta^q\le\sup_{r\ge0}\left\{\frac{\beta g}Wr-\frac{\kappa_q}{q} r^q\right\}
   =\frac{q-1}{q}\kappa_q^{-1/(q-1)}
       \left(\frac{\beta g}W\right)^{q/(q-1)}.
\end{equation}
Consequently $g\ge A_*h^{(q-1)/q}$. 
During a good iteration \eqref{s6:stepconditions} we have
\begin{equation}
 h_k-h_{k+1}\ge \varsigma \vartheta g_k^{(1+\alpha)/\alpha} \ge \varsigma \vartheta A_*^{(1+\alpha)/\alpha}h^{(q-1)(1+\alpha)/(q\alpha)} \ge c_*h_k^{1+s}.
 \label{s6:recurrence}
\end{equation}
For non-good iterations, we have that $h_{k+1}\le h_k$.
If $s>0$ and $h_{k+1}>0$, then
\begin{equation*}
 h_{k+1}^{-s}-h_k^{-s} =s\int_{h_{k+1}}^{h_k}r^{-s-1}\,dr \ge s(h_k-h_{k+1})h_k^{-s-1}\ge sc_*
\end{equation*}
at each good iteration. 
Summing over all steps proves the convergence rate for $s>0$. 
For $s=0$, \eqref{s6:recurrence} gives $h_{k+1}\le(1-\theta)h_k$ at good iterations, which is the standard iteration bound for linear convergence rate.
\end{proof}

\begin{remark}
\label{s6:badsteps}
As shown in \cite[Appendix~D]{lacoste2015global}, before termination, any block of consecutive iterations that are not good has length at most $B_N:=3N!$. In particular,
\begin{equation}
 t_k\ge\left\lfloor\frac{k}{B_N+1}\right\rfloor.
 \label{s6:count}
\end{equation}
\end{remark}

To conclude, we show that the convergence result holds when the objective function is the composition between a convex function and a linear function.
This extends the quadratic conditioning argument in \cite[Appendix~F]{lacoste2015global} to arbitrary powers $q\ge2$.

\begin{proposition}
\label{s6:hoffman}
Suppose $F(x)=\mathcal H(Ax)+\langle b,x\rangle$ on $P$ and that $\mathcal H$ is continuously differentiable on a neighborhood of $A(\Isubl)$.
Assume that $\mathcal H$ is $q$-uniformly convex on $A(\Isubl)$ with constant $\mu_H>0$ and $q\ge2$, as in \eqref{eq:growth_cond_intro}.
Then $F$ is $q$-conditioned on $\Isubl$ as in \eqref{eq:pcond_intro}, and there exists a suitable constant $\kappa_q$ for which \eqref{eq:unif_p_cond_assumption} holds with $\beta=2$.    
\end{proposition}

\begin{proof}
The strict convexity of $\mathcal H$ implies that all optimizers have the same image $y^\star=Ax^\star$.
Moreover, following the argument proposed by \cite{beck2017linearly}, we have $b^\star=\langle b,x^\star\rangle$.
These two equalities characterize optimality among points of $P$, hence
\begin{equation*}
 \XC^\star=\{x\in P:Ax=y^\star,\ \langle b,x\rangle=b^\star\}.
\end{equation*}
Hoffman's error bound for this consistent system of linear equalities and the inequalities defining $P$ gives a finite constant $K>0$ such that
\begin{equation*}
 \operatorname{dist}(x,\XC^\star)
 \le K\big(\|Ax-y^\star\|+|\langle b,x\rangle-b^\star|\big)
 \qquad(x\in P).
\end{equation*}
Let $L_H=\max_{y\in A(C)}\|\nabla\mathcal H(y)\|$, $\bar h=F(x_0)-F^\star$, $u=\|Ax-y^\star\|$, and $h=F(x)-F^\star$. 
Convexity of $C$ ensures that $\mathcal H$ is $L_H$-Lipschitz on $A(C)$, whence
\begin{equation*}
 |\langle b,x\rangle-b^\star|\le h+L_Hu,
 \qquad
 \operatorname{dist}(x,\XC^\star)^q
 \le C_H\big((1+L_H)^qu^q+\bar h^{q-1}h\big),
 \quad C_H:=2^{q-1}K^q.
\end{equation*}
At an optimizer, first-order optimality and \eqref{eq:growth_cond_intro} give $h\ge(\mu_H/q)u^q$. 
Hence \eqref{eq:pcond_intro} holds with
\begin{equation}
 \nu:=\left[C_H\left((1+L_H)^q\frac q{\mu_H}
                                  +\bar h^{q-1}\right)\right]^{-1}.
 \label{s6:nu}
\end{equation}
For the reverse Bregman divergence, put
$a_x=-\langle\nabla F(x),x^\star-x\rangle$. The linear part cancels and
\begin{equation*}
 a_x-h=\mathcal H(y^\star)-\mathcal H(Ax)
             -\langle\nabla\mathcal H(Ax),y^\star-Ax\rangle
 \ge\frac{\mu_H}qu^q.
\end{equation*}
Therefore $-h+2a_x\ge h+(2\mu_H/q)u^q$. Combining the last two estimates
proves \eqref{eq:unif_p_cond_assumption} with $\beta=2$ and, for example,
\begin{equation}
 \kappa_q:=\frac{q}{C_H\max\{\bar h^{q-1},\,
                                 q(1+L_H)^q/(2\mu_H)\}}>0.
 \label{s6:kappa}
\end{equation}
\end{proof}

\subsection{Support Identification Property}
\label{sec:support_identification}

In this section, we extend the support identification property presented in \cite{bomze2020active} to fit our framework.
Since the proof techniques are close to the one used in \cite{bomze2020active}, we only report the main result that will be used in what follows and defer the proof in a dedicated part of the Appendix.
First, we introduce some useful notation.
For each $i\in [n]$, we define $\theta_i(x)=\scalarp{\nabla F(x), e_i - x}$. 
Given $x^\star\in\mathcal{X}^\star$, we denote by $I(x^\star)$ the extended support of $x^\star$, i.e., the set of indices $i\in [n]$ for which $\theta_i(x^\star)=0$. 
By definition, the complement of $I(x^\star)$, that is $I^c(x^\star)=[n]\setminus I(x^\star)$, corresponds to the set of indices satisfying $\theta_i(x^\star)>0$. 
Given $x_0$ and fixed a stepsize, let $x^{(k)}$ be the $k$-th iteration of the PFW algorithm, then set
\[
O_k(x^\star) := \{i\in I^c(x^\star) \;|\; x_i^{(k)}=0\} \qquad \text{and}\qquad J_k(x^\star)=I^c(x^\star)\setminus O_k(x^\star).
\]
Finally, we define $\delta_k=\max_{[n]\setminus O_k(x^\star)}\theta_i(x^\star)=\max_{I(x^\star)\cup J_k(x^\star)}\theta_i(x^\star)$.
When the element $x^\star\in\mathcal{X}^\star$ is clear from the context or is unique, we simplify the notation by omitting its dependence, e.g., $I=I(x^\star)$, $I^c=I^c(x^\star)$, and so on.
For the sake of simplicity, we present the result for the case in which $P$ is the $n$ dimensional simplex and extend the results to the general case afterwards.

\begin{theorem}
\label{thm:identification_properties}
    Let $F$ be a continuously differentiable function whose gradient is $\alpha$-Hölder continuous with constant $H$ over an $\varepsilon$-neighborhood of the set of optimal solutions in $\Delta_n$.
    For every $k$, let $x^{(k)}$ and $d^{(k)}$ be the $k$-th iterate and the direction found at $x^{(k)}$ by the PFW algorithm.
    Denoted by $\gamma_k$ the step size at step $k$, assume that for every $k$ either 
    \begin{equation*}
        \gamma_k=\gamma^{\rm max}_k\qquad\text{or}\qquad\gamma_k\ge \Big(- \frac{\scalarp{\nabla f(x^{(k)}),d^{(k)}}}{2H}\Big)^{\frac{1}{\alpha}}.
    \end{equation*}
    Finally, given $x^\star\in\XC^\star$, set
    \begin{equation*}
        \delta_{\rm min}:=\min \big\{\theta_i(x^\star)\;\text{s.t.}\;\theta_i(x^\star)>0\big\}\quad\text{and}\quad\delta_{\mathrm{max}}:=\max \big\{\theta_i(x^\star)\;\text{s.t.}\;\theta_i(x^\star)>0\big\}.
    \end{equation*}
    Then, if $\|x^{(k)}-x^\star\|_1 < r^\star=\min\Big\{\varepsilon,2^{\frac{\alpha-1}{\alpha}}\Big(\frac{\delta_{\rm min}}{2Hn^{\frac{1}{2}}+\delta_{\mathrm{max}}}\Big)^{\frac{1}{\alpha}}\Big\}$, we have
    \begin{equation*}
        |J_{k+1}|\le \max\big\{|J_k|-1,0\big\}.
    \end{equation*}
\end{theorem}

\begin{remark}[Extending the results to the general polytope]
    The results presented in Theorem \ref{thm:identification_properties} and Lemma \ref{lmm:ext_Hölder_lambda} extend to the case in which $\mathcal{P}$ is a general polytope as follows.
    Given the finite set of vertices of $\mathcal{P}$, namely $\{v^{(i)}\}_{i\in[M]}$ take any linear map $A_\mathcal{P}:\Delta_M\to\mathcal{P}$ for which it holds 
    \begin{equation}
        A_\mathcal{P}(w)= \sum_{i=1}^Mw_iv^{(i)}
    \end{equation}
    where $w\in\Delta_M$.
    Since minimizing $F$ over $\mathcal{P}$ is equivalent to minimizing $F_{A_\mathcal{P}} = F \circ A_\mathcal{P}$ over $\Delta_M$, Theorem~\ref{thm:identification_properties} applies directly, thereby extending the identification result to the general setting in which $\mathcal{P}$ is an arbitrary polytope.
\end{remark}

\subsection{Application to UOT problems induced by entropy like divergences}
In this subsection, we apply the results obtained for the PFW algorithm to the UOT problems induced by $\phi=\Up$.
We recall that $\phipz=\Upo$, $\phipi=\Upi$, and
\begin{equation*}
 \Up^{\prime}(t)=\frac{t^{p-1}-1}{p-1}\qquad\text{and}\qquad \Up^{\prime\prime}(t)=t^{p-2}\qquad(t\in\R_{++}).
\end{equation*}
Via a classic calculus argument, we have the following.

\begin{lemma}
\label{s6:scalar}
Let $M>0$ and set $\DivBp(y,x)=\Up(y)-\Up(x)-\Up'(x)(y-x)$, then
\begin{enumerate}[label=(\roman*)]
\item If $p\ge2$ for $x,y\in\left[0,M\right]$, it holds
\begin{equation*}
 \DivBp(y,x)\ge\frac{|y-x|^p}{p(p-1)},\qquad
 |\Up^\prime(y)-\Up'(x)|\le M^{p-2}|y-x|.
\end{equation*}
\item If $1<p\le2$ for $x,y\in\left[0,M\right]$, it holds
\begin{equation*}
 \DivBp(y,x)\ge\frac{M^{p-2}}2|y-x|^2,\qquad
 |\Up^\prime(y)-\Up'(x)|\le\frac{|y-x|^{p-1}}{p-1}.
\end{equation*}
\item If $p\le1$ and $0<m_0<M$, the derivative is Lipschitz and the
function is strongly convex on $\left[m_0,M\right]$, with respective
constants $m_0^{p-2}$ and $M^{p-2}$.
\end{enumerate}
\end{lemma}

\begin{proof}
For $p\ge2$ and $y\ge x$, writing $d=y-x$ gives
\begin{equation*}
 \DivBp(y,x)=\int_0^d(d-t)(x+t)^{p-2}\,dt
 \ge\int_0^d(d-t)t^{p-2}\,dt=\frac{d^p}{p(p-1)}.
\end{equation*}
A similar argument handles the case in which $y<x$, concluding the first part of the statement.
For $1<p\le2$, we have $\Up^{\prime\prime}(t)\ge M^{p-2}$ on $\left]0,M\right]$, which proves the quadratic bound by integration. 
The inequality $|a^r-b^r|\le|a-b|^r$ for $a,b\ge0$ and $0<r\le1$ follows from $(a+b)^r\le a^r+b^r$, and proves the derivative estimate with $r=p-1$.
The last statement of the theorem follows by integrating $M^{p-2}\le\Up^{\prime\prime}(t)\le m_0^{p-2}$.
\end{proof}

For $p>1$, we then have that $\Up$ is both $q=\max\{2,p\}$ conditioned and has $\min\{p-1,1\}$-Hölder continuous gradient on any set of the form $[0,M]$.
If $p\le0$, the derivative of $\psi_p$ is not $\alpha$-Hölder on any set of the form $]0,M]$.
However, since $\Up(t)\to+\infty$ as $t\downarrow0$, every finite initial sub-level set is bounded away from zero at every positive-reference marginal allowing us to bypass problematic points by considering a suitable sub-level set.
The remaining case $0<p\le1$ the sub-level set approach does not apply since $\Up$ has bounded values in $0$ but its derivative is not $\alpha$-Hölder on any finite sub-level set. 
To bypass this limitation, we introduce an additive barrier term to $\Up$ that preserves the solutions to the induced UOT problem and allows us to handle the resulting entropies as done for $p\le 0$.

\begin{proposition}
\label{s6:barrier}
Let $R>0$ and $0<p\le1$, then there exists a positive value $t_0$ such that every optimizer of $\GR$ on $\setKR$ satisfies $r_i/\mu_i\ge t_0$ if $\mu_i>0$ and $s_j/\nu_j\ge t_0$ if $\nu_j>0$.
Then, if $\delta=t_0/2$ and $\zeta>0$, we set
\begin{equation}
 b_\delta(t)=
 \begin{cases}\log(\delta/t)+t/\delta-1,&0<t<\delta,\\
 0,&t\ge\delta,\\+\infty,&t\le0,
 \end{cases}
 \qquad \widetilde{\Up}=\Up+\zeta b_\delta.
 \label{s6:amendment}
\end{equation}
Replacing every positive-reference perspective induced by $\Up$ with the one induced by $\widetilde{\Up}$, we obtain an objective $\widetilde{\GR}\ge \GR$ with the same minimum and exactly the same minimizers on $\setKR$. 
Moreover $\widetilde{\Up}$ has a Lipschitz gradient and is strongly convex on each finite sub-level set.
\end{proposition}

\begin{proof}
The first part of the proof follows from Theorem \ref{thm:lower-marginals} with $U_i=V_j=2B_0$.
Let us then consider the function $b_\delta$ defined as in \eqref{s6:amendment}.
By definition, $b_\delta$ is nonnegative, convex and $C^1$ on
$\left]0,+\infty\right[$.
If $t<\delta$ the derivative of $b_\delta$ is $-1/t+1/\delta$ and its second derivative is $1/t^2$, so that $b_\delta$ is continuous and has a continuous first derivative.
Since $b_\delta\ge 0$ we have $\widetilde{\Up}\ge \Up$, thus $\widetilde{\GR}\ge \GR\ge v_R$ for every $(\eta,u,w)\in\FeasR$.
Since $b_\delta(t)=0$ for every $t\ge \delta$, we have that $\widetilde{\GR}(\eta^\star,u^\star,w^\star) = \GR(\eta^\star,u^\star,w^\star) = v_R$ for every $(\eta^\star,u^\star,w^\star)\in\argmin \GR$.
We then infer that $(\eta^\star,u^\star,w^\star)\in\argmin \widetilde{\GR}$.
Conversely, a minimizer of $\widetilde{\GR}$ must also minimize $\GR$, since $\GR\le\widetilde{\GR}$ and $\min\widetilde{\GR}=v_R$.

To conclude, notice that on any finite sub-level set, every positive-reference marginal is bounded away from zero, while all marginals are at most $2B_0$.
On such compact intervals the derivative of $\widetilde{\Up}$ is Lipschitz and its curvature is bounded below by the positive
minimum of $\Up^{\prime\prime}$.
Scaling the perspectives and summing proves the last assertion, concluding the proof.
\end{proof}

Proposition \ref{s6:barrier} allows us to apply the results presented for the generic PFW algorithm to every UOT problem induced by an entropy-like penalty function.

\begin{theorem}
\label{thm:convergence_SPFW}
Given $P=\setKR$, for $p\in\R$ consider
\begin{equation}
 F_p=\begin{cases}
 \GR,\quad\quad&p\le0\text{ or }p>1,\\
 \widetilde{\GR},&0<p\le1
 \end{cases}
 \qquad\qquad
 \text{and}
 \qquad\qquad
 F_p^\star=v_{p,R}.
 \label{s6:workingobjective}
\end{equation}
Let $x_0$ be a finite-value initial point and let $x^{(k)}$ be the sequence generated by the PFW algorithm combined with a stepsize as in Proposition~\ref{s6:linesearch}.
Then, the hypotheses of Theorems~\ref{thm:general_convergence_PFW} and \ref{thm:identification_properties} hold, and the following rates apply
\begin{equation}
    \label{s6:polyrate}
    F_p(x^{(k)})-v_{p,R} =
\begin{cases}
O\left(t_k^{-\frac{p}{p-2}}\right),
    & p > 2, \\[6pt]
O\left((1-\theta)^{t_k}\right),
    & p = 2, \\[6pt]
O\left(t_k^{-\frac{2(p-1)}{2-p}}\right),
    & 1 < p < 2, \\[6pt]
O\left((1-\theta_p)^{t_k}\right),
    & p \leq 1,
\end{cases}
\end{equation}
where $\theta_p\in(0,1)$ depends on $p\le 1$ and $t_k$ is the number of good steps made up to iterate $k$.
\end{theorem}

\begin{proof}
Let $\mathsf{A}$ be the linear function defined as
\[
    \mathsf{A}(z)=\mathsf{A}(\eta,u,w)=(\Pruno\eta+u,\Prdue\eta+w).\footnote{If necessary, restrict the vector to include only indices associated with positive entries.}
\]
Then, we have $F(z)=\mathcal H(\mathsf{A}z)+\langle b,z\rangle$, for a suitable vector $b$ and $\mathcal H$ is the sum of perspective functions associated with a positive $\mu_i$ or a positive $\nu_j$ and induced by a suitable $q$-uniform convex function.
From our assumptions, $\mathcal{H}$ is not a null function since at least one of $\mu$ and $\nu$ is not null and, for a given reference weight $m$, the perspective is $\varphi_m(t)=m\Up(t/m)$. 
Moreover, for each $\varphi_m$ with $m>0$, we have $\varphi_m^{\prime\prime}(t)=m^{1-p}t^{p-2}$ so that Lemma~\ref{s6:scalar} ensures $\varphi_m$ is $q=\max\{2,p\}$ conditioned over the set $[0,2B_0]$ if $p>1$.
For $p\le 1$, we can restrict our attention to a sub-level set to ensure the objective is $2$-conditioned.
On the same sets, the function has also Lipschitz continuous gradient if $p\ge 2$, Hölder continuous if $1 < p< 2$ or Lipschitz continuous on the sub-level set.

These properties readily translate to $\mathcal{H}$ since it is the sum of the perspectives.
Indeed, denoted by $d_+>0$ the number of positive entries, we have
\begin{equation*}
 \|\nabla\mathcal H(y)-\nabla\mathcal H(z)\|^2
 \le L_0^2\sum_{i=1}^{d_+}|y_i-z_i|^{2\alpha}
 \le L_0^2d_+^{1-\alpha}\|y-z\|^{2\alpha},
\end{equation*}
where $\alpha\in]0,1]$ and $L_0$ is a suitable constant that depends on $\Up$.
In particular, $F$ has Hölder continuous gradient with positive constant $L_0d_+^{(1-\alpha)/2}\|A\|^{1+\alpha}$.
Also, since $\mathcal{H}$ is the sum of independent $q$-uniform convex functions, it is also $q$-uniform convex.
Thus Proposition \ref{s6:hoffman} allows us to conclude that $F$ satisfies Assumption \ref{assumption_F_P} and thus Theorem \ref{thm:general_convergence_PFW} is applicable.
The convergence rates follow by substituting the appropriate values $(q,\alpha)$ in \eqref{s6:rateconstants}. 
\end{proof}

\begin{remark}[Generality beyond the power family]
While we presented our results for powerlike entropies, the same approach used by us can be applied to any UOT problem in which the function $\phi$ induces a Csiszár divergence that satisfies Assumption \ref{assumption_F_P} and whose composed gradient is $\alpha$-Hölder continuous. 
\end{remark}

\subsection{Reconstruction and a computable stopping certificate}

In this section, we show how to use the reconstruction functions $\Lift,\Bsf,\Dsf$ to attain an optimality certificate that allow to relate the minimum of $\GR$ and $F_c$.

\begin{proposition}[Certificates for inexact compressed iterates]
\label{s6:reconstruct}
Let $z=(\eta,u,w)\in\setKR$ have finite objective value. 
For each $\Xsf\in\{\Lift,\Dsf,\Bsf\}$, we have
\begin{equation}
 F_c(\Xsf(z))-\widetilde{\GR}(z)\le\mathcal R_{\Xsf}(z).
 \label{s6:recoverydifference}
\end{equation}
Moreover, 
\begin{equation}
 0\le F_c(\Xsf(z))-v_c
 \le \widetilde{\GR}(z)-v_R+\mathcal R_{\Xsf}(z)
 \le g_R^{\rm FW}(z)+\mathcal R_{\Xsf}(z),
 \label{s6:certificate}
\end{equation}
where $v_R$ the minimum objective of $\widetilde{\GR}$ on $\setKR$ and $g_R^{\rm FW}(z)$ is the Frank-Wolfe gap of the objective $\widetilde{\GR}$ at point $z$.
\end{proposition}

\begin{proof}
The inequality in \eqref{s6:recoverydifference} follows from the same argument used in Theorem \ref{thm:error-certificates} applied to a generic element $(\eta,u,w)\in\FeasR$.
Then, since $v_R\le v_c$, we have
\begin{equation*}
 F_c(\Xsf(z))-v_c\le\widetilde{\GR}(z)-v_R+\mathcal R_{\Xsf}(z)
                       \le g_R^{\rm FW}(z)+\mathcal R_{\Xsf}(z),
\end{equation*}
where in the last inequality we used the standard Frank-Wolfe gap bound.
\end{proof}

Notice that a similar result holds for $\widetilde{\GR}$.
In the end, we are able to evaluate the available certificates and return the corresponding reconstructed plan as soon as
\begin{equation}
 g_R^{\rm FW}(z_k)+\min_{\Xsf\in\{\Lift,\Bsf,\Dsf\}}\mathcal R_{\Xsf}(z_k)
 \le\varepsilon.
 \label{s6:stop}
\end{equation}
This certifies the objective error for the original objective $F_c$ via quantities that depend on $R$ and/or on $\widetilde{\GR}$.
To conclude, we use these optimality certificates to provide a bound on the number of iterates needed to identify the support of the solution.

\begin{theorem}
\label{thm:complexity_lambda_general_p}
Suppose the hypotheses of Theorem~\ref{thm:general_convergence_PFW} hold for a working objective $\widetilde{\GR}$ whose minimum is $v_R$. 
Define
\begin{equation}
 E_{B_0}(R):=B_0\min\left\{
 \begin{array}{ll}
 {[c_{\max}-R]_+}, &\\
 {[-2\phipz-R]_+}, &\text{if }\phipz> -\infty,\\
 {[\max_i\minr+\max_j\minc+2\phipi-R]_+},
                   &\text{if }\phipi<+\infty.
 \end{array}\right\},
 \label{s6:uniformbudget}
\end{equation}
take $\varepsilon>0$, $0<\varrho\le1$ such that $E_{B_0}(R)\le(1-\varrho)\varepsilon$, and any starting point $x_0$. 
If $c_0\ge h_0> \varrho\varepsilon$, we have that after $T(\varrho\varepsilon)$ good steps, reconstruction by a method attaining the minimum in \eqref{s6:uniformbudget} gives $F_c(\pi)-v_c\le\varepsilon$, where
\begin{equation}
 T(\varrho\varepsilon):=
 \begin{cases}
 \displaystyle\left\lceil\frac{(\varrho\varepsilon)^{-s}-c_0^{-s}}{sc_*}\right\rceil,
        &s>0,\\[6pt]
 \displaystyle\left\lceil\frac{\log(c_0/\varrho\varepsilon)}{-\log(1-\theta)}\right\rceil,
        &s=0.
 \end{cases}
 \label{s6:complexity}
\end{equation}
\end{theorem}

\begin{proof}
Since $\tau\le B_0$, Proposition~\ref{s6:reconstruct} bounds each error by a term in \eqref{s6:uniformbudget}.
For the redirecting plan use $\sum_i\minr u_i\le\tau\max_i\minr$ and $\sum_j\minc w_j\le\tau\max_j\minc$. 
For the selected reconstruction, \eqref{s6:certificate} gives $F_c(\pi)-v_c\le h_k+E_{B_0}(R)$.

If $s>0$, $h_0\le c_0$ and \eqref{s6:polyrate} imply $h_k\le(c_0^{-s}+sc_*t_k)^{-1/s}$. 
Solving for this bound to be at most $\varrho\varepsilon$ gives the first line of \eqref{s6:complexity}. 
Solving $c_0(1-\theta)^{t_k}\le \varrho\varepsilon$ gives the second line. 
The total error is then at most $\varrho\varepsilon+(1-\varrho)\varepsilon=\varepsilon$.
Remark~\ref{s6:badsteps} proves the actual-iteration count.
Early termination means $h_k=0$ and hence also suffices.
\end{proof}

\printbibliography


\appendix

\section{Proof of the support identification property}

\begin{lemma}
\label{lmm:ext_Hölder_lambda}
    Given a continuously differentiable function $F$ over the simplex $\Delta_n$, let $x^\star\in\Delta_n$ be a minimizer of $F$.
    Given $h>0$, let $x^{(k)} \in \Delta_{n}$ be such that $\|x^\star-x^{(k)}\|_1\le h$, then, if $F$ has $\alpha$-Hölder continuous gradient over $B_{h,1}(x^\star)\cap\Delta_n$, we have
    \begin{equation}
    \label{eq:hsimple}
        |\theta_i(x^\star)-\theta_i(x^{(k)})|\le h^\alpha\Big(Hn^{\frac{1}{2}}+\frac{\delta_kh^{1-\alpha}}{2}\Big),
    \end{equation}
    where $H$ is the Hölder constant associated with the gradient and $\delta_k=\max_{i\in[n]\backslash O_k}\theta_i(x^\star)$ 
\end{lemma}

\begin{remark}
    Notice that, if $\alpha=1$, the result is the same as the one presented in \cite{bomze2020active}.
    Moreover, if $h<1$, equation \eqref{eq:hsimple} simplifies into
    \begin{equation}
        |\theta_i(x^\star)-\theta_i(x^{(k)})|\le h^\alpha\Big(Hn^{\frac{1}{2}}+\frac{\delta_k}{2}\Big),
    \end{equation}
    which is equivalent to the condition presented in \cite{bomze2020active}, whose Lipschitz constant is computed with respect to the $\|\cdot\|_1$ norm.
\end{remark}

\begin{proof}
    By definition of $\theta_i(x)$, we have
    \begin{align*}
        |\theta_i(x^\star)-\theta_i(x^{(k)})|&=|\scalarp{\nabla F(x^{(k)})-\nabla F(x^\star), e_i-x^{(k)}} -\scalarp{\nabla F(x^\star), x^{(k)}-x^\star}|\\
        &\le |\scalarp{\nabla F(x^{(k)})-\nabla F(x^\star), e_i-x^{(k)}}| + |\scalarp{\nabla F(x^\star), x^{(k)}-x^\star}|\\
        &\le \|\nabla F(x^{(k)})-\nabla F(x^\star)\|\|e_i-x^{(k)}\| + |\scalarp{\nabla F(x^\star), x^{(k)}-x^\star}|\\
        &\le H \|x^\star-x^{(k)}\|^\alpha n^{\frac{1}{2}} + |\scalarp{\nabla F(x^\star), x^{(k)}-x^\star}|\\
        &\le H n^{\frac{1}{2}} h^\alpha + |\scalarp{\nabla F(x^\star), x^{(k)}-x^\star}|,
    \end{align*}
    where in the last inequality we have used the hypothesis on $\|x^\star-x^{(k)}\|_1$.
    We now provide an estimate of the term $|\scalarp{\nabla F(x^\star), x^{(k)}-x^\star}|$.
    Since both $x^{(k)}$ and $x^\star$ belong to $\Delta_n$, we have
    \begin{equation}
        0=\sum_{i=1}^n(x^{(k)}_i-x^\star_i) \quad \iff \quad \sum_{i=1}^n(x^{(k)}_i-x^\star_i)_+=\sum_{i=1}^n(x^{(k)}_i-x^\star_i)_-,   
    \end{equation}
    where $(t)_+$ and $(t)_-$ denote the positive and negative part of $t$.
    In particular, for every $\ell\in[n]$, 
    \begin{equation}
        |x^{(k)}_\ell-x^\star_\ell| \le \sum_{i=1}^n(x^{(k)}_i-x^\star_i)_+=\frac{1}{2}\|x^{(k)}-x^\star\|_1\le \frac{h}{2}.
    \end{equation}
    Finally, since $x^\star$ is optimal, we have $|\scalarp{\nabla F(x^\star),x^{(k)}-x^\star}|\le (\nabla F(x^\star)_M-\nabla F(x^\star)_m)\frac{h}{2}\le \delta_k\frac{h}{2}$, where $\nabla F(x^\star)_M$ and $\nabla F(x^\star)_m$ are the maximum and minimum entries of $\nabla F$ over $[n]\backslash O_k$, respectively.
    We then conclude
    \begin{equation}
    \label{eq:simplehno1-alpha}
        |\theta_i(x^\star)-\theta_i(x^{(k)})|\le Hn^{\frac{1}{2}}h^\alpha+\frac{\delta_k}{2}h= h^\alpha\Big(Hn^{\frac{1}{2}}+\frac{\delta_kh^{1-\alpha}}{2}\Big).
    \end{equation}
\end{proof}

\begin{theorem}
    Let $F$ be a continuously differentiable function whose gradient is $\alpha$-Hölder continuous with constant $H$ over an $\varepsilon$-neighborhood of the set of optimal solutions in $\Delta_n$.
    For every $k$, let $x^{(k)}$ and $d^{(k)}$ be the $k$-th iterate and the direction found at $x^{(k)}$ by the PFW algorithm.
    Denoted by $\gamma_k$ the step size at step $k$, assume that for every $k$ either 
    \begin{equation*}
        \gamma_k=\gamma^{\rm max}_k\qquad\text{or}\qquad\gamma_k\ge \Big(- \frac{\scalarp{\nabla f(x^{(k)}),d^{(k)}}}{2H}\Big)^{\frac{1}{\alpha}}.
    \end{equation*}
    Finally, given $x^\star\in\XC^\star$, set
    \begin{equation*}
        \delta_{\rm min}:=\min \big\{\theta_i(x^\star)\;\text{s.t.}\;\theta_i(x^\star)>0\big\}\quad\text{and}\quad\delta_{\mathrm{max}}:=\max \big\{\theta_i(x^\star)\;\text{s.t.}\;\theta_i(x^\star)>0\big\}.
    \end{equation*}
    Then, if $\|x^{(k)}-x^\star\|_1 < r^\star=\min\Big\{\varepsilon,2^{\frac{\alpha-1}{\alpha}}\Big(\frac{\delta_{\rm min}}{2Hn^{\frac{1}{2}}+\delta_{\mathrm{max}}}\Big)^{\frac{1}{\alpha}}\Big\}$, we have
    \begin{equation*}
        |J_{k+1}|\le \max\big\{|J_k|-1,0\big\}.
    \end{equation*}
\end{theorem}

\begin{proof}
    Let us denote by $J_k$ the set $\{i\in[n]\;\text{s.t.}\;\theta_i(x^\star)>0\;\text{and}\;x^{(k)}_i>0\}$.    
    If $I^c=\emptyset$, then $\theta_i(x^\star)=0$ for every $i\in[n]$, thus $J_k=\emptyset$ for every $k$.
    Let us then assume that $I^c\neq\emptyset$. 
    First, we notice that by definition $\frac{\delta_{\rm min}}{2Hn^{\frac{1}{2}}+\delta_{\mathrm{max}}}\le 1$, hence \eqref{eq:simplehno1-alpha} holds.
    In this case $\delta_{\mathrm{max}}\ge\delta_k\ge \delta_{\rm min}>0$.
    Owing to Lemma \ref{lmm:ext_Hölder_lambda}, for every $i\in I^c$ it holds
    \begin{align*}
        \theta_i(x^{(k)}) &\ge \theta_i(x^\star)-\|x^{(k)}-x^\star\|_1^\alpha\Big(Hn^{\frac{1}{2}}+\frac{\delta_k}{2}\Big) > \theta_i(x^\star)-\frac{\delta_{\rm min}}{2Hn^{\frac{1}{2}}+\delta_{\mathrm{max}}}\Big(Hn^{\frac{1}{2}}+\frac{\delta_k}{2}\Big)\\
        &=
        \theta_i(x^\star)-\frac{\delta_{\rm min}}{2}\frac{2Hn^{\frac{1}{2}}+\delta_k}{2Hn^{\frac{1}{2}}+\delta_{\mathrm{max}}},
    \end{align*}
    since $2^{\frac{\alpha-1}{\alpha}}\le 1$ for every $\alpha\in ]0,1]$.
    We now split the study into two cases, depending on whether $J_k=\emptyset$ or $J_k\neq \emptyset$.
    \textbf{Case 1, $J_k=\emptyset$.} 
    Since $J_k=\emptyset$, we have that $O_k=I^c$, thus $x^{(k)}_i=0$ for every $i\in I^c$ and hence $\delta_k=0$.
    In this case, for every $i\in I^c$
    \begin{equation}
    \label{eq:lambda_pos_alt}
        \theta_i(x^{(k)}) \ge \theta_i(x^\star)-\frac{Hn^{\frac{1}{2}}\delta_{\rm min}}{2Hn^{\frac{1}{2}}+\delta_{\mathrm{max}}}\ge \theta_i(x^\star)-\frac{Hn^{\frac{1}{2}}\delta_{\rm min}}{2Hn^{\frac{1}{2}}}\ge \delta_{\rm min}-\frac{\delta_{\rm min}}{2}>0.
    \end{equation}

    By definition of the pairwise direction, we then have that any $i\in I^c$ cannot be an away step direction as $i$ does not belong to the active set of the solution $x^{(k)}$.
    Similarly, since $\theta_i(x^{(k)})>0$ (see \eqref{eq:lambda_pos_alt}), none of the indices $i\in I^c$ can be selected as a standard Frank-Wolfe direction.
    Hence $x^{(k+1)}_i=0$ for every $i\in I^c$.

    \textbf{Case 2, $J_k\neq\emptyset$.} 
    By definition of $\delta_k$, we have that there exists an index $\ell$ for which 
    \[
        \theta_\ell(x^\star)=\max_{j\in J_k} \theta_j(x^\star)=\delta_k.
    \]
    Owing again to Lemma \ref{lmm:ext_Hölder_lambda}
    \begin{align*}
        \theta_\ell(x^{(k)}) &\ge \theta_\ell(x^\star)-\|x^{(k)}-x^\star\|_1^\alpha\Big(Hn^{\frac{1}{2}}+\frac{\delta_k}{2}\Big) > \delta_{k}-\frac{\delta_{\rm min}}{2Hn^{\frac{1}{2}}+\delta_{\mathrm{max}}}\Big(Hn^{\frac{1}{2}}+\frac{\delta_k}{2}\Big)\\
        &\ge \delta_{k}-\frac{\delta_{\rm min}}{2}\frac{2Hn^{\frac{1}{2}}+\delta_k}{2Hn^{\frac{1}{2}}+\delta_{\rm min}}
        \ge \delta_{k}-\frac{\delta_{k}}{2}\frac{2Hn^{\frac{1}{2}}+\delta_k}{2Hn^{\frac{1}{2}}+\delta_{k}} = \frac{\delta_{k}}{2},
    \end{align*}
    where we have used the fact that the function $a\to-\frac{a}{2Hn^{1/2}+a}$ is decreasing on $[0,+\infty[$ and $\delta_k\ge \delta_{\rm min}$. 
    Notice that, as a by-product, we get
    \begin{equation}
        \label{eq:byproduct}
        \frac{\delta_{\rm min}}{2}\frac{2Hn^{\frac{1}{2}}+\delta_k}{2Hn^{\frac{1}{2}}+\delta_{\rm min}}\le \frac{\delta_k}{2}.
    \end{equation}

    We now show that the pairwise direction selected at step $k$, namely $d^{(k)}=e_i-e_j$, is such that $i\in I$ and $j\in J_k$.
    First, we notice that for every $i\in I^c$, we have $\theta_i(x^{(k)})>0$, thus the FW direction must belong to $I$.
    Second, we show that $j\in J_k$ (or, equivalently, that $j\notin I$ since the indices of $O_k$ do not belong to the active set of the solution $x^{(k)}$).
    For every $j\in I$, we have that $\theta_j(x^\star)=0$, thus
    \begin{align}
        |\theta_j(x^{(k)})|&=|\theta_j(x^{(k)})-\theta_j(x^\star)|\le \|x^{(k)}-x^\star\|_1^\alpha\Big(Hn^{\frac{1}{2}}+\frac{\delta_k}{2}\Big)\\
        &<\frac{\delta_{\rm min}}{2Hn^{\frac{1}{2}}+\delta_{\mathrm{max}}}\Big(Hn^{\frac{1}{2}}+\frac{\delta_k}{2}\Big)\le\frac{\delta_{\rm min}}{2Hn^{\frac{1}{2}}+\delta_{\rm min}}\Big(Hn^{\frac{1}{2}}+\frac{\delta_k}{2}\Big)\le \frac{\delta_k}{2},
    \end{align}
    where in the last inequality, we have used \eqref{eq:byproduct}.
    We then have that none of the directions in $I$ can determine the away step direction, since $\frac{\delta_k}{2}<\theta_\ell(x^{(k)})$.
    To conclude, we need to show that at each iteration the algorithm takes a full step, that is $\gamma_k=\gamma^{(\mathrm{max})}_k$.
    Toward a contradiction, let us assume that $\gamma_k<\gamma^{(\mathrm{max})}_k$.
    By hypothesis,
    \[
        \gamma_k\ge \Big(-\frac{\scalarp{\nabla F(x^{(k)}),e_i-e_j}}{2H}\Big)^{\frac{1}{\alpha}} =
        \Big(\frac{\theta_j(x^{(k)})- \theta_i(x^{(k)}) }{2H}\Big)^{\frac{1}{\alpha}}
    \]
    owing to the fact that $\theta_i(x^{(k)})\le 0$, we infer
    \[
        \gamma_k\ge\Big(\frac{\theta_j(x^{(k)})}{2H}\Big)^{\frac{1}{\alpha}}\ge\Big(\frac{\delta_{\rm min}}{4H}\Big)^{\frac{1}{\alpha}}.
    \]
    Since $j\in I^c$, we have $x_j^\star=0$, thus
    \[
        x_j^{(k)}=(x^{(k)}-x^\star)_j\le \frac{\|x^{(k)}-x^\star\|_1}{2}\le\frac{1}{2^{\frac{1}{\alpha}}}\Big(\frac{\delta_{\rm min}}{2Hn^{\frac{1}{2}}+\delta_{\mathrm{max}}}\Big)^{\frac{1}{\alpha}}=\Big(\frac{\delta_{\rm min}}{4Hn^{\frac{1}{2}}+2\delta_{\mathrm{max}}}\Big)^{\frac{1}{\alpha}}.
    \]

    Therefore, by the algorithm definition
    \begin{equation}
        x_j^{(k+1)}=x^{(k)}_j-\gamma_k\le \Big(\frac{\delta_{\min}}{4Hn^{\frac{1}{2}}+2\delta_{\mathrm{max}}}\Big)^{\frac{1}{\alpha}}-\Big(\frac{\delta_{\rm min}}{4H}\Big)^{\frac{1}{\alpha}}<0,
    \end{equation}
    which is a contradiction.
\end{proof}

From Theorem \ref{thm:identification_properties} we infer that, whenever we use a PFW algorithm to solve an optimization problem, the routine identifies the extended support of the solution in a finite number of iterations.
The number of iterations depends on two terms: \begin{enumerate*}[label=(\roman*)]
    \item the number of iterations needed to satisfy the condition $\|x^{(k)}-x^\star\|_1\le r^\star$, where $x^\star$ is an optimal point, and
    \item the number of steps needed to reduce the cardinality of the set $J_{k+1}:=\{i\in[n]\;\text{s.t.}\;\theta_i(x^\star)>0\;\text{and}\;x^{(k+1)}_i>0\}$ to zero, which is at most equal to the dimension of the simplex. 
\end{enumerate*}

\clearpage

\end{document}